\documentclass[a4paper,oneside]{article}
\usepackage[utf8]{inputenc}
\usepackage[a4paper,margin=3cm]{geometry} %showframe=true
\usepackage{amsmath}
\usepackage{amsthm}
\usepackage{amscd}
\usepackage{amssymb}
\usepackage{MnSymbol}
\usepackage{latexsym}
\usepackage{eucal}
\usepackage{dsfont}
\usepackage{mathtools}
\usepackage{enumitem}
\usepackage{todonotes}
\usepackage{verbatim}
\usepackage{graphicx}
\usepackage{float}
\usepackage{array,booktabs}

\usepackage{pdflscape}

\usepackage[colorlinks,pdftex]{hyperref}
\hypersetup{
linkcolor=black,
citecolor=black,
pdftitle={},
pdfauthor={Franziska K\"uhn},
pdfkeywords={},
}

\makeatletter
\renewcommand\section{\@startsection{section}{1}{0mm}{-1.5\baselineskip}{\baselineskip}{\normalsize\bfseries\sffamily}}
\renewcommand\subsection{\@startsection{subsection}{1}{0mm}{-\baselineskip}{\baselineskip}{\normalsize\bfseries\sffamily}}
\makeatother

\makeatletter
\def\@fnsymbol#1{\ensuremath{\ifcase#1\or *\or **\or \dagger\or \ddagger\or
   \mathsection\or \mathparagraph\or \|\or \dagger\dagger
   \or \ddagger\ddagger \else\@ctrerr\fi}}

\newlength{\preskip}
\newlength{\postskip}
\makeatother

\newtheoremstyle{theorem}{\preskip}{\postskip}{\itshape}{}{\bfseries}{}
{.5em}{\textbf{\thmname{#1}\thmnumber{ #2} (\thmnote{ #3})}}
\newtheoremstyle{definition}{\preskip}{\postskip}{\normalfont}{0pt}{\bfseries}{}{.5em}{}
\newtheoremstyle{remark}{\preskip}{\postskip}{\normalfont}{0pt}{\bfseries}{}{.5em}{}

\swapnumbers
\theoremstyle{theorem} \newtheorem{thm}{Theorem}[section]
\theoremstyle{theorem} 
\theoremstyle{theorem} 
\theoremstyle{theorem} \newtheorem{kor}[thm]{Corollary}
\theoremstyle{definition} 
\theoremstyle{remark} 
\theoremstyle{remark} 
\theoremstyle{definition} 
\theoremstyle{definition} \newtheorem*{ack}{Acknowledgements}
\theoremstyle{remark} 
\theoremstyle{remark} \newtheorem{bems}[thm]{Remarks}
\theoremstyle{definition}  \newtheorem{bsp}[thm]{Example}
\theoremstyle{definition}

\DeclareMathOperator \re {Re}
\DeclareMathOperator \im {Im}

\DeclareMathOperator \tr {tr}

\newcommand{\I}{\mathds{1}}

\newcommand\fa{\qquad \text{for all \ }}

\newcommand{\cadlag}{c\`adl\`ag }

\newcommand\mc[1] {\mathcal{#1}}
\newcommand\mbb[1] {\mathds{#1}}

\newcommand\T{\rule{0pt}{3.5ex}}       % Top strut
\newcommand\B{\rule[-3ex]{0pt}{0pt}} % Bottom strut

\author{%
    Franziska K\"{u}hn\thanks{Institut f\"ur Mathematische Stochastik, Fachrichtung Mathematik, Technische Universit\"at Dresden, 01062 Dresden, Germany, \texttt{franziska.kuehn1@tu-dresden.de}} 
}

\title{Transition probabilities of L\'evy-type processes: Parametrix construction}

\date{}

\begin{document}

\maketitle

\abstract{\noindent We present an existence result for L\'evy-type processes which requires only weak regularity assumptions on the symbol $q(x,\xi)$ with respect to the space variable $x$.  Applications range from existence and uniqueness results for L\'evy-driven SDEs with H\"{o}lder continuous coefficients to existence results for stable-like processes and L\'evy-type processes with symbols of variable order. Moreover, we obtain heat kernel estimates for a class of L\'evy and L\'evy-type processes. The paper includes an extensive list of L\'evy(-type) processes satisfying the assumptions of our results. \par \medskip

  \par \medskip

}

\section{Introduction} \label{intro}

The L\'evy--Khintchine formula gives a one-to-one correspondence between L\'evy processes and continuous negative definite functions (i.\,e.\ characteristic exponents). For any continuous negative definite function $\psi: \mbb{R}^d \to \mbb{C}$ with $\psi(0)=0$, i.\,e.\ any function of the form \begin{equation*}
	\psi(\xi) =\psi(0) -ib \cdot \xi + \frac{1}{2} \xi \cdot Q \xi + \int_{y \neq 0} (1-e^{i\xi \cdot y} +i \xi \cdot y \I_{(0,1)}(|y|)) \, \nu(dy), \qquad \xi \in \mbb{R}^d
\end{equation*}
for some $b \in \mbb{R}^d$, a positive semidefinite symmetric matrix $Q \in \mbb{R}^{d \times d}$ and a measure $\nu$ on $\mbb{R}^d \backslash \{0\}$ such that $\int_{\mbb{R}^d \backslash \{0\}} |y|^2 \wedge 1 \, \nu(dy)<\infty$, there exists a L\'evy process with characteristic exponent $\psi$. Over the past years, there has been an increasing interest in so-called L\'evy-type processes. This is a class of Markov processes which behave locally like a L\'evy process, but the (analogue of the) L\'evy triplet depends on the current position of the process in the state space. If the smooth functions with compact support $C_c^{\infty}(\mbb{R}^d)$ are contained in the domain of the generator of a L\'evy-type process $(X_t)_{t \geq 0}$, then the process can be characterized via its symbol $q$, \begin{equation}
	q(x,\xi) =  q(x,0)-ib(x) \cdot \xi + \frac{1}{2} \xi \cdot Q(x) \xi + \int_{y \neq 0} (1-e^{i \xi \cdot \xi}+i \xi \cdot y \I_{(0,1)}(|y|)) \, \nu(x,dy), \quad x,\xi \in \mbb{R}^d, \label{intro-eq1}
\end{equation}
which is a continuous negative definite function for each fixed $x \in \mbb{R}^d$. Typical examples are processes with variable index of stability (this corresponds to $q(x,\xi) = |\xi|^{\alpha(x)}$) and solutions of L\'evy-driven SDEs, see Table~\ref{tab-ltp} on page \pageref{tab-ltp} for further examples. \par
It is natural to ask whether for a given function $q$ of the form \eqref{intro-eq1}, i.\,e.\ family $(q(x,\xi))_{x \in \mbb{R}^d}$ of continuous negative definite functions, there exists a L\'evy-type process with symbol $q$. The answer is, in general, \emph{no} (see e.\,g.\ \cite[Example 2.26]{ltp} for counterexamples), and therefore it is important to find sufficient conditions on the symbol $q$ or the characteristics $(b(x),Q(x),\nu(x,dy))_{x \in \mbb{R}^d}$ which ensure the existence of a L\'evy-type process with a given symbol $q$. Many existence results in the literature (see \cite{ltp} for an overview) are rather restrictive in the sense that they either assume that $q$ is of a particular form (typically ``stable-like'' or ``lower-order perturbation'') or they require strong assumptions on the regularity of the symbol $q$ with respect to the space variable $x$ (typically smoothness). \par
In this paper, we present a new existence result for L\'evy-type processes which requires only mild regularity assumptions on the symbol $q$ with respect to $x$. Applications range from variable order subordination and L\'evy-type processes with symbols of variable order to existence and uniqueness results for solutions to L\'evy-driven stochastic differential equations (SDEs) with globally H\"{o}lder continuous coefficients. The proof of the main result, Theorem~\ref{ltp-3}, relies on a parametrix construction; this is an analytical method which has become an increasingly popular tool in probability theory to prove the existence of certain stochastic processes and to establish heat kernel estimates,  e.\,g.\  processes with variable order of differentiation (Kolokoltsov \cite{kol00,kol} and Chen \& Zhang \cite{zhang}), gradient perturbations of L\'evy generators (Bogdan \& Jakubowski \cite{bogdan07} and Jakubowski \& Szczypkowski \cite{jak09}) and solutions of SDEs with H\"{o}lder continuous coefficients (Knopova \& Kulik \cite{kul15,kul15-2} and Huang \cite{huang15}). Using the parametrix construction, we can not only prove the existence of L\'evy-type processes with symbols from certain classes, but also get additional information on the process such as heat kernel estimates, well-posedness of the associated martingale problem or the richness of the domain of the generator, cf.\ Theorem~\ref{ltp-5}. \par
As a by-product of the parametrix construction, we obtain heat kernel estimates for a class of L\'evy processes; the estimates are crucial to prove the convergence of the parametrix expansion. Estimates for the transition density of L\'evy processes have attracted a lot of attention, for example heat kernel estimates for unimodal \cite{bogdan13}, rotationally invariant \cite{chen13}, tempered stable \cite{szt10} L\'evy processes or L\'evy processes with exponential moments \cite{knop12,szt15}, to mention but a few. In contrast to many of these results, we state our assumptions on the L\'evy process in terms of the characteristic exponent $\psi$ and not in terms of the L\'evy triplet $(b,Q,\nu)$. Since there are few Lévy processes for which both the characteristic exponent and the L\'evy triplet can be calculated explicitly, both approaches (i.\,e.\ via L\'evy triplet or via characteristic exponent) have their own justification. In fact, for most of the examples which we present in Section~\ref{levy} the L\'evy triplet is unknown, and therefore it is very hard to verify conditions on the L\'evy triplet. Our result applies, for instance, to relativistic stable, Lamperti stable, normal tempered stable, truncated L\'evy processes, cf.\ Example~\ref{levy-7}, and a class of subordinators; see Table~\ref{tab-lp} and Table~\ref{tab-sub} in Section~\ref{levy} for further examples. \par \medskip

This paper is organized as follows. In Section~\ref{def} we introduce the basic definitions and some notation. Section~\ref{levy} is devoted to heat kernel estimates for L\'evy processes. In Section~\ref{ltp} we present our main result, the existence result for L\'evy-type processes, and discuss several applications, including a new existence and uniqueness result for L\'evy-driven SDEs with H\"{o}lder continuous coefficients (Corollary~\ref{ltp-15}), an existence result for stable-like processes (Example~\ref{ltp-11}) and Feller processes with symbols of variable order (Corollary~\ref{ltp-17}). \par
The results presented in this paper are essentially taken from my PhD thesis \cite{diss}. The aim of this paper is to give a brief summary of the most important results and make them accessible to a larger audience; in particular, we do not include proofs since they are very technical and lengthy, and we refer to \cite{diss,matters} for full proofs.

\section{Preliminaries} \label{def}

We consider the Euclidean space $\mbb{R}^d$ endowed with the canonical scalar product $x \cdot y = \sum_{j=1}^d x_j y_j$ and the Borel-$\sigma$-algebra $\mc{B}(\mbb{R}^d)$. The continuous bounded functions are denoted by $C_b(\mbb{R}^d)$, and $C_{\infty}(\mbb{R}^d)$ is the space of continuous functions $f: \mbb{R}^d \to \mbb{R}$ vanishing at infinity. Superscripts $k\in\mbb{N}$ are used to denote the order of differentiability, e.\,g.\ $f \in C_{\infty}^k(\mbb{R}^d)$ means that $f$ and its derivatives up to order $k$ are $C_{\infty}(\mbb{R}^d)$-functions. A function $\ell: (0,\infty) \to (0,\infty)$ is \emph{slowly varying (at infinity)} if \begin{equation*}
	\lim_{x \to \infty} \frac{\ell(\lambda x)}{\ell(x)}=1 \fa \lambda>0.
\end{equation*}
We say that a function $f: \mbb{R}^d \to \mbb{R}^k$ is (globally) H\"{o}lder continuous if there exist constants $\varrho \in (0,1]$, $C>0$ such that $|f(x)-f(y)| \leq C |x-y|^{\varrho}$ for all $x,y \in \mbb{R}^d$. \par
A $d$-dimensional Markov process $(\Omega,\mc{A},\mbb{P}^x,x \in \mbb{R}^d,X_t,t \geq 0)$ with \cadlag (right-continuous with left-hand limits) sample paths is called a \emph{L\'evy-type process} if the associated semigroup $(P_t)_{t \geq 0}$ defined by \begin{equation*}
	P_t f(x) := \mbb{E}^x f(X_t), \quad x \in \mbb{R}^d, \; f \in \mc{B}_b(\mbb{R}^d) := \{f: \mbb{R}^d \to \mbb{R}; \text{$f$ bounded, Borel measurable}\}
\end{equation*}
has the \emph{Feller property} and $(P_t)_{t \geq 0}$ is \emph{strongly continuous at $t=0$}, i.\,e. $P_t f \in C_{\infty}(\mbb{R}^d)$ for all $C_{\infty}(\mbb{R}^d)$ and $\|P_tf-f\|_{\infty} \xrightarrow[]{t \to 0} 0$ for any $f \in C_{\infty}(\mbb{R}^d)$. L\'evy-type processes are also known as \emph{Feller processes}; we use both terms synonymously. A semigroup $(P_t)_{t \geq 0}$ has the \emph{strong Feller property} if $P_t f \in C_b(\mbb{R}^d)$ for any $f \in \mc{B}_b(\mbb{R}^d)$, $t \geq 0$. If the smooth functions with compact support $C_c^{\infty}(\mbb{R}^d)$ are contained in the domain of the generator $(L,\mc{D}(L))$, then we speak of a \emph{rich} L\'evy-type process. A result due to von Waldenfels and Courr\`ege, cf.\ \cite[Theorem 2.21]{ltp}, states that the generator $L$ of a rich L\'evy-type process is, when restricted to $C_c^{\infty}(\mbb{R}^d)$, a pseudo-differential operator with negative definite symbol: \begin{equation*}
	Lf(x) =  - \int_{\mbb{R}^d} e^{i \, x \cdot \xi} q(x,\xi) \hat{f}(\xi) \, d\xi, \qquad f \in C_c^{\infty}(\mbb{R}^d), \, x \in \mbb{R}^d,
\end{equation*}
where $\hat{f}(\xi) := (2\pi)^{-d} \int_{\mbb{R}^d} e^{-ix \xi} f(x) \, dx$ denotes the Fourier transform of $f$ and \begin{equation}
	q(x,\xi) = q(x,0) - i b(x) \cdot \xi + \frac{1}{2} \xi \cdot Q(x) \xi + \int_{\mbb{R}^d \backslash \{0\}} \left(1-e^{i y \cdot \xi}+ i y\cdot \xi \I_{(0,1)}(|y|)\right) \, \nu(x,dy). \label{cndf}
\end{equation}
We call $q$ the \emph{symbol} of the rich L\'evy-type process $(X_t)_{t \geq 0}$ and $-q$ the symbol of the pseudo-differential operator. For each fixed $x \in \mbb{R}^d$, $(b(x),Q(x),\nu(x,dy))$ is a L\'evy triplet, i.\,e.\ $b(x) \in \mbb{R}^d$, $Q(x) \in \mbb{R}^{d \times d}$ is a symmetric positive semidefinite matrix and $\nu(x,dy)$ a $\sigma$-finite measure on $(\mbb{R}^d \backslash \{0\},\mc{B}(\mbb{R}^d \backslash \{0\}))$ satisfying $\int_{y \neq 0} \min\{|y|^2,1\} \, \nu(x,dy)<\infty$. A set $\mc{D} \subseteq C_{\infty}(\mbb{R}^d)$ is a \emph{core} for the generator $(L,\mc{D}(L))$ if $\overline{(L,\mc{D})}^{\|\cdot\|_{\infty}} = (L,\mc{D}(L))$. Our standard reference for L\'evy-type processes is the monograph \cite{ltp}. \par
A \emph{L\'evy process} $(L_t)_{t \geq 0}$ is a rich Feller process whose symbol $q$ does not depend on $x$. This is equivalent to saying that $(L_t)_{t \geq 0}$ has stationary and independent increments and \cadlag sample paths, cf.\ \cite[Theorem 2.6]{ltp}. The symbol $q=q(\xi)$ (also called \emph{characteristic exponent}) and the L\'evy process $(L_t)_{t \geq 0}$ are related through the L\'evy--Khintchine formula: \begin{equation*}
 \mbb{E}^xe^{i \xi \cdot (L_t-x)} = e^{-t q(\xi)} \fa t \geq 0, \, x,\xi \in \mbb{R}^d.
\end{equation*}
A L\'evy process with non-decreasing sample paths is a \emph{subordinator} and can be characterized by its Laplace exponent, cf.\ \cite{bernstein}. We refer to Sato \cite{sato} for a detailed discussion of L\'evy processes and to Schilling \cite{barca} for an introduction to L\'evy and L\'evy-type processes.

\section{Heat kernel estimates for L\'evy processes} \label{levy}

In this section, we present transition density estimates for a class of L\'evy processes. The main results are Theorem~\ref{levy-1} (heat kernel estimates for rotationally invariant L\'evy processes in dimension $d \geq 1$), Theorem~\ref{levy-3} (heat kernel estimates for one-dimensional L\'evy processes which are not necessarily symmetric) and Corollary~\ref{levy-9} (heat kernel estimates for subordinators).  In contrast to many results in the literature, we state our assumptions on the L\'evy process $(L_t)_{t \geq 0}$ in terms of the characteristic exponent and not in terms of the L\'evy triplet. 

\begin{thm} \label{levy-1}
	Let $(L_t)_{t \geq 0}$ be a $d$-dimensional L\'evy process and suppose that its characteristic exponent $\psi: \mbb{R}^d \to \mbb{C}$ satisfies \eqref{L1}-\eqref{L3}.
\begin{enumerate}[label*=\upshape (L\arabic*),ref=\upshape L\arabic*] 
			\item\label{L1} $\psi$ is rotationally invariant, i.\,e.\ there exists $\Psi: \mbb{R} \to \mbb{R}$ such that $\psi(\xi) = \Psi(|\xi|)$, $\xi \in \mbb{R}^d$. If $m>0$: $\Psi$ is even, i.\,e.\ $\Psi(r) = \Psi(-r)$ for all $r \geq 0$.
			\item\label{L2} There exists $\theta \in (0,\tfrac{\pi}{2})$ and $m \geq 0$ such that $\Psi$ has a holomorphic extension to \begin{equation}
				\Omega := \Omega(m,\theta) := \{z \in \mbb{C}; |\im z|<m\} \cup \{z \in \mbb{C} \backslash \{0\}; \arg z \in (-\theta,\theta) \cup (\pi-\theta,\pi+\theta)\}. \label{domain}
			\end{equation}
					\begin{figure}[H]
						\begin{center}
							\includegraphics[scale=0.75]{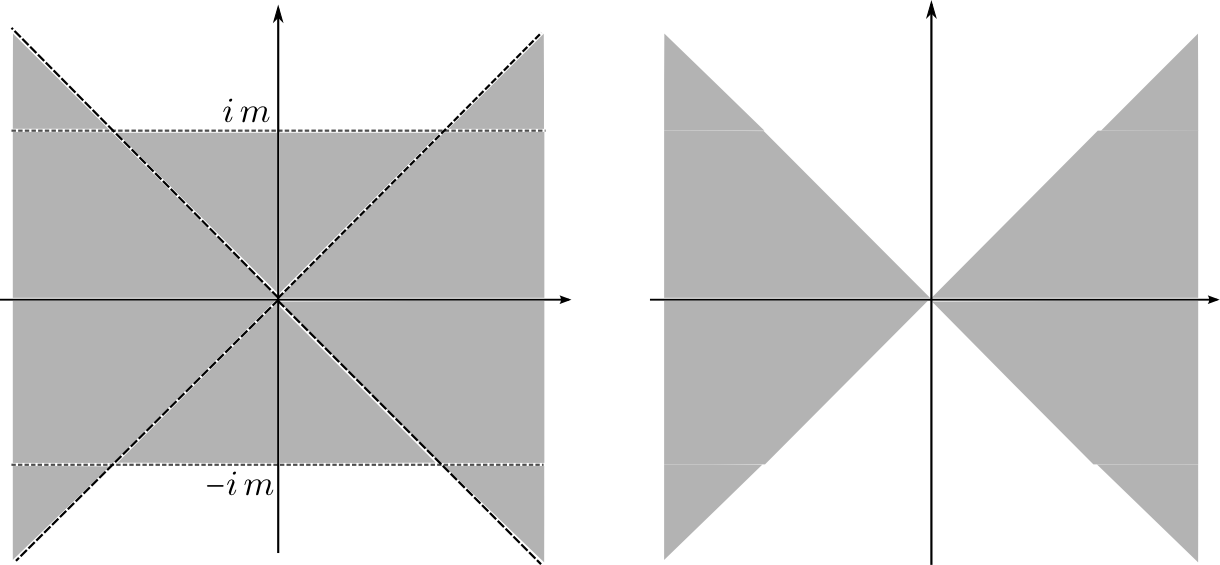}
						\end{center}
						\caption{The domain $\Omega = \Omega(m,\vartheta)$ for $m>0$ (left) and $m=0$ (right).}
						\label{fig:def_gebiet_exp}
					\end{figure} 
		\item\label{L3} There exist constants $c_1,c_2>0$ and $\gamma_0,\gamma_{\infty} \in (0,2]$ and a slowly varying (at infinity) increasing function $\ell: (0,\infty) \to (0,\infty)$ such that \begin{equation}
				\re \Psi(z) \geq ´\frac{c_1}{\ell(|z|)} |z|^{\gamma_{\infty}} \fa z \in \Omega, |z| \geq 1 \label{l3-eq1} 
			\end{equation}
		and \begin{equation}
				|\Psi(z)| \leq c_2 \ell(|z|) \left(|z|^{\gamma_0} \I_{\{|z| \leq 1\}} + |z|^{\gamma_{\infty}} \I_{\{|z| >1\}}\right) \fa z \in \Omega. \label{l3-eq2} 
			\end{equation}
	\end{enumerate}
	Then the transition density \begin{equation*}
			p_t(x) := \frac{1}{(2\pi)^d} \int_{\mbb{R}^d} e^{-ix \cdot \xi} e^{-t \psi(\xi)} \, d\xi
	\end{equation*}
	exists, is infinitely often differentiable and satisfies the estimates \begin{align*}
			|p_t(x)| &\leq C S(x,t) \\
			\left| \frac{\partial}{\partial t} p_t(x) \right| &\leq C t^{-1} S(x,t)
		\end{align*}
	for any $x \in \mbb{R}^d$, $t \in (0,T]$ where $C=C(T)>0$ is an absolute constant and \begin{equation*}
			S(x,t) := S_m(x,t) := \exp \left( - \frac{m}{4} |x| \right) (1+\ell(ct^{-1/\gamma_{\infty}})) \begin{cases} t^{-d/\gamma_{\infty}}, & |x| \leq t^{1/\gamma_{\infty}} \wedge 1, \\ t/|x|^{d+\gamma_{\infty}}, & t^{1/\gamma_{\infty}} < |x| \leq 1, \\ t/|x|^{d+\gamma_{\infty} \wedge \gamma_0}, & |x|>1\end{cases}
		\end{equation*}
		for some absolute constant $c=c(T)>0$. 	Moreover, for any multi-index $\beta=(\beta_1,\ldots,\beta_d) \in \mbb{N}_0^d$ and any $T>0$ there exists a constant $c>0$ such that \begin{equation*}
			\left| \frac{\partial^{\beta}}{\partial x^{\beta}} p_t(x) \right| \leq c t^{-(|\beta_1|+\ldots+|\beta_d|)/\gamma_{\infty}} S(x,t) \fa x \in \mbb{R}^d, t \in (0,T]. 
	\end{equation*}
\end{thm}

In dimension $d=1$ we do not need any symmetry assumption like rotational invariance in \eqref{L1}.

\begin{thm} \label{levy-3}
	Let $(L_t)_{t \geq 0}$ be a one-dimensional L\'evy process with characteristic exponent $\psi$. If $\Psi(\xi) := \psi(\xi)$ satisfies \eqref{L2} and \eqref{L3}, then the results of Theorem~\ref{levy-1} remain valid.
\end{thm}

To prove Theorem~\ref{levy-1} and Theorem~\ref{levy-3} we apply Cauchy's theorem to shift the contour of integration; for details see \cite[Section 4.1]{matters} or \cite[Section 4.1]{diss}. Using this idea, we can also obtain heat kernel estimates for subordinators, cf.\ Corollary~\ref{levy-9}. Let us give some remarks on both Theorem~\ref{levy-1} and Theorem~\ref{levy-3}.

\begin{bems} \label{levy-5}\begin{enumerate}
	\item Condition \eqref{L3} implies that $\psi$ satisfies the \emph{sector condition}, i.\,e.\ there exists a constant $C_1>0$ such that \begin{equation*}
		|\im \psi(\xi)| \leq C_1 |\re \psi(\xi)| \fa \xi \in \mbb{R}^d.
	\end{equation*}
	Moreover, it follows from \eqref{L3} that there exists a constant $C_2>0$ such that \begin{equation*}
		\re \Psi(z) \geq - C_2>-\infty \fa z \in \Omega.
	\end{equation*}
	Note that the existence of such a constant is not trivial; although any continuous negative definite function $\psi: \mbb{R}^d \to \mbb{C}$ satisfies $\re \psi(\xi) \geq 0$ for all $\xi \in \mbb{R}^d$, the inequality does, in general, not need to be true for \emph{complex} $z$.
	\item The constant $m/4$ in the definition of the function $S$ introduced in Theorem~\ref{levy-1} can be replaced by $m(1-\delta)$ for any $\delta \in (0,1)$. Since it is well-known, cf.\ \cite[Theorem 5.26]{ushakov}, that the characteristic function of $L_t$ \begin{equation*}
		\mbb{E}e^{i \xi L_t} = e^{-t \psi(\xi)}, \qquad t > 0, \xi \in \mbb{R}^d,
	\end{equation*}
	is analytic on the strip $\{z \in \mbb{C}; |\im z|<m\}$ if, and only if, $\mbb{E}e^{m(1-\delta)|L_t|}<\infty$ for any $\delta \in (0,1)$, the exponential decay $\exp(-|x| m(1-\delta))$ is, in general, the best we can expect.
	\item Let $\psi$ be a continuous negative definite function with L\'evy triplet $(b,0,\nu)$. If $\psi$ satisfies the sector condition, then it is possible to give sufficient conditions in terms of fractional moments of $\nu|_{B(0,1)}$ and $\nu|_{B(0,1)^c}$ which ensure that $\psi$ satisfies the growth condition \eqref{l3-eq2} for \emph{real} $z$, cf.\ Blumenthal \& Getoor \cite{blumen} and Schilling \cite{rs97}. This is, however, no longer possible for the holomorphic extension. For instance if we consider $\psi(\xi) := 1-\cos \xi$, then the associated L\'evy measure $\nu = \tfrac{1}{2} \delta_1 + \frac{1}{2} \delta_{-1}$ has arbitrary moments, but the  (unique) holomorphic extension $\Psi(z) = 1-\cos z$ does not satisfy \eqref{l3-eq2}. 
\end{enumerate} \end{bems}

If the characteristic exponent $\psi$ is given in closed form, it is usually easy to check whether the assumptions of Theorem~\ref{levy-1} are satisfied.

\begin{bsp} \label{levy-7}
	Let $(L_t)_{t \geq 0}$ be a  $d$-dimensional L\'evy process with one of the following characteristic exponents. \begin{enumerate}
			\item\label{levy-7-i} (isotropic stable) $\psi(\xi) = |\xi|^{\alpha}$, $\xi \in \mbb{R}^d$, $\alpha \in (0,2]$,
			\item\label{levy-7-ii} (relativistic stable) $\psi(\xi) = (|\xi|^2+\varrho^2)^{\alpha/2}-\varrho^{\alpha}$, $\xi \in \mbb{R}^d$, $\varrho>0$, $ \alpha \in (0,2)$,
			\item\label{levy-7-iii} (Lamperti stable) $\psi(\xi) = (|\xi|^2+\varrho)_{\alpha}-(\varrho)_{\alpha}$, $\xi \in \mbb{R}^d$, $\varrho>0$, $\alpha \in (0,1)$, where $(r)_{\alpha} := \Gamma(r+\alpha)/\Gamma(r)$ denotes the Pochhammer symbol,
			\item\label{levy-7-iv} (TLP: truncated L\'evy process) $\psi(\xi) = (|\xi|^2+\varrho^2)^{\alpha/2} \cos(\alpha \arctan(\varrho^{-1} |\xi|))-\varrho^{\alpha}$, $\xi \in \mbb{R}^d$, $\alpha \in (0,2)$, $\varrho>0$,
			\item\label{levy-7-v} $\psi(\xi) = (|\xi|^{\beta}-1)/(|\xi|^{\alpha}-1)$, $\xi \in \mbb{R}^d$, $0<\alpha <\beta \leq 1$, 
			\item\label{levy-7-vi} (NTS: normal tempered stable, $d=1$) $\psi(\xi) = (\kappa^2+(\xi-ib)^2)^{\alpha/2}-(\kappa^2-b^2)^{\alpha/2}$, $\xi \in \mbb{R}$, $\alpha \in (0,2)$, $b>0$, $\kappa>b$.
		\end{enumerate}
		The characteristic exponents \ref{levy-7-i}-\ref{levy-7-v} and \ref{levy-7-vi} satisfy the assumptions of Theorem~\ref{levy-1} and Theorem~\ref{levy-3}, respectively, with \begin{enumerate}
				\item $\gamma_0 = \gamma_{\infty} = \alpha$, $m=0$, $\ell=1$, 
				\item $\gamma_0 = 2$, $\gamma_{\infty} = \alpha$, $m = \varrho$, $\ell=1$, 
				\item $\gamma_0 = 2$, $\gamma_{\infty} = 2 \alpha$, $m = \sqrt{\varrho}$, $\ell=1$,
				\item $\gamma_0 = 2$, $\gamma_{\infty} = \alpha$, $m = \varrho$, $\ell=1$,
				\item $\gamma_0 = \alpha$, $\gamma_{\infty} = \beta-\alpha$, $m=0$, $\ell=1$,
				\item $\gamma_0 = 2$, $\gamma_{\infty} = \alpha$, $m = \kappa-b$, $\ell=1$.
			\end{enumerate}
		Consequently, Theorem~\ref{levy-1}, respectively, Theorem~\ref{levy-3}, provide heat kernel estimates for the transition density $p$, its time derivative and derivatives with respect to the space variable $x$. For the particular case of isotropic stable L\'evy processes, we recover well-known (sharp) estimates for the heat kernel, cf.\ \cite{bendikov,blumen60}. %ref rel. stable
	\end{bsp}
				
The following list of examples satisfying the assumptions of Theorem~\ref{levy-1} and Theorem~\ref{levy-3}, respectively, is taken from \cite[Table 5.2]{matters}.

\begin{landscape} 
\begin{table} \centering
	\begin{tabular}{llllll}
		& \textbf{Name} & \textbf{char.\ Exponent} & \textbf{Parameter} & \textbf{Dim.} & \textbf{Heat Kernel Estimate (Thm.~\ref{levy-1})}  \\
	\hline
		1
			& isotropic $\alpha$-stable
			& $\displaystyle|\xi|^{\alpha}$  
			& $\alpha \in (0,2]$
			& $d \geq 1$
			&$\gamma_0=\gamma_{\infty}=\alpha$, $m=0$, $\ell=1$ \B\T \\
		\hline
		2
			& relativistic stable
			& $\displaystyle(|\xi|^2+\varrho^2)^{\alpha/2}-\varrho^{\alpha}$
			& $\alpha \in (0,2)$, $\varrho \in (0,\infty)$
			& $d \geq 1$
			&$\gamma_0=2$, $\gamma_{\infty}=\alpha$, $m \in (0,\varrho)$, $\ell=1$ \B\T \\ \hline
3
	& \vtop{ \hbox{\strut normal tempered stable} \hbox{\strut (NTS)}}
	& $\displaystyle(\kappa^2+(\xi-ib)^2)^{\alpha/2}-(\kappa^2-b^2)^{\alpha/2}$
	& \vtop{ \hbox{\strut$\alpha \in (0,2)$, $\kappa \in (0,\infty)$,} \hbox{\strut$b<\kappa$}}
	& $d=1$
	& $\gamma_0=2$, $\gamma_{\infty}=\alpha$, $m \in (0,\kappa-b)$, $\ell=1$ \B\T \\ \hline
4
	& 
	& $\displaystyle\frac{|\xi|^2}{\sqrt{|\xi|^2+\varrho}}$
	& $\varrho \in (0,\infty)$
	& $d \geq 1$
	& $\gamma_0=2$, $\gamma_{\infty}=1$, $m \in (0,\sqrt{\varrho)}$, $\ell=1$ \B\T \\ \hline
5
	& 
	& $\displaystyle\frac{|\xi|^2}{(|\xi|^2+\varrho)^{\alpha}}$
	& $\alpha \in (0,2)$, $\varrho \in (0,\infty)$
	& $d \geq 1$
	& $\gamma_0=2$, $\gamma_{\infty}=2-\alpha$, $m \in (0,\sqrt{\varrho})$, $\ell=1$\B\T \\ \hline
6
	& 
	& $\displaystyle\frac{|\xi|^{\beta}-1}{|\xi|^{\alpha}-1}-1$ (extended by continuity at $\xi=1$)
	& $0<\alpha<\beta < 1$
	& $d \geq 1$
	& $\gamma_0=\alpha$, $\gamma_{\infty}=\beta-\alpha$, $m=0$, $\ell=1$\B\T \\ \hline
7
	& 
	& $\displaystyle-\frac{|\xi|^{\alpha}-1}{|\xi|^{\alpha-2}-1}$ (extended by continuity at $\xi=1$)
	& $\alpha \in (0,2)$
	& $d \geq 1$
	& $\gamma_0=2-\alpha$, $\gamma_{\infty}=2$, $m =0$, $\ell=1$ \B\T \\ \hline
8
	& 
	& $\displaystyle\frac{|\xi|^{\alpha}-1}{|\xi|^{\alpha-2}-1}-1$ (extended by continuity at $\xi=1$)
	& $\alpha \in (2,4]$
	& $d \geq 1$
	& $\gamma_0=\alpha-2$, $\gamma_{\infty}=2$, $m=0$, $\ell=1$ \B\T \\ \hline
9
	& 
	& $\displaystyle|\xi|^2 \frac{|\xi|^{\alpha}-\varrho^{\alpha}}{|\xi|^2-\varrho^2}$ (extended by continuity at $\xi=\varrho$)
	&$\alpha \in (0,2)$, $\varrho \in (0,\infty)$
	& $d \geq 1$
	& $\gamma_0=2$, $\gamma_{\infty}=\alpha$, $m=0$, $\ell=1$\B\T \\ \hline
10
	& 
	& $\displaystyle(|\xi|^{-\alpha}+|\xi|^{-\beta})^{-1}$ (extended by continuity at $\xi=0$)
	& $\alpha,\beta \in (0,2]$
	& $d \geq 1$
	& $\gamma_0=\alpha \vee \beta$, $\gamma_{\infty}=\alpha \wedge \beta$, $m=0$, $\ell=1$ \B\T \\ \hline
11
	& 
	& $\displaystyle|\xi| (1-e^{-2 \varrho |\xi|})$
	& $\varrho \in (0,\infty)$
	& $d \geq 1$
	& $\gamma_0=2$, $\gamma_{\infty}=1$, $m=0$, $\ell=1$ \B\T \\ \hline
12
	& 
	& $\displaystyle|\xi| (1+e^{-2 \varrho |\xi|})$
	& $\varrho \in (0,\infty)$
	& $d \geq 1$
	& $\gamma_0=\gamma_{\infty}=1$, $m=0$, $\ell=1$ \B\T \\ \hline
13
	& 
	& $\displaystyle\varrho |\xi|^2 (|\xi|^2+1) \log(1+|\xi|^{-2})$ %|\xi|^2 \log(1+|\xi|^{-2}) \to 1$ as $|\xi| \to \infty
	& $\varrho \in (0,\infty)$
	& $d \geq 1$
	& $\gamma_0=\gamma_{\infty}=2$, $m=0$, $\ell(r)=\log(r \vee e)$ \end{tabular} \caption{Examples of L\'evy processes satisfying the assumptions of Theorem~\ref{levy-1} and Theorem~\ref{levy-3}, respectively} \label{tab-lp}. \end{table}
	\setcounter{table}{0}
\begin{table} \centering
\begin{tabular}{llllll}
	& \textbf{Name} & \textbf{char.\ Exponent} & \textbf{Parameter} & \textbf{Dim.} & \textbf{Heat Kernel Estimate (Thm.~\ref{levy-1})}  \\ \hline
14
	& 
	& $\displaystyle \varrho \frac{|\xi|^2 (|\xi|^2+1)}{(|\xi|^2+2) \log(|\xi|^2+2)}$
	& $\varrho \in (0,\infty)$
	& $d \geq 1$
	& $\gamma_0=\gamma_{\infty}=2$, $m \in (0,\sqrt{2})$, $\ell(r)=\log(r \vee e)$ \B\T \\ \hline
15
	& 
	& $\displaystyle |\xi| \arctan(\varrho |\xi|)$
	& $\varrho \in (0,\infty)$
	& $d \geq 1$
	& $\gamma_0=\gamma_{\infty}=1$, $m=0$, $\ell=1$ \B\T \\ \hline
16
	& \vtop{ \hbox{\strut truncated L\'evy process} \hbox{\strut (TLP)}}
	& $\displaystyle (|\xi|^2+\varrho^2)^{\alpha/2} \cos \big(\alpha \arctan \tfrac{|\xi|}{\varrho}\big)-\varrho^{\alpha}$
	& $\alpha \in (0,2)$, $\varrho \in (0,\infty)$
	& $d \geq 1$
	& $\gamma_0=2$, $\gamma_{\infty}=\alpha$, $m \in (0,\varrho)$, $\ell=1$ \B\T \\ \hline
17
	& 
	& $\displaystyle \varrho |\xi| \frac{\cosh^2(\sqrt{2} |\xi|)}{\sinh(2 \sqrt{2} |\xi|)}$
	& $\varrho \in (0,\infty)$
	& $d \geq 1$
	& $\gamma_0=\gamma_{\infty}=1$, $m=0$, $\ell=1$ \B\T \\ \hline
18
	& 
	& $\displaystyle \varrho |\xi| \frac{\sinh^2(\sqrt{2} |\xi|)}{\sinh(2 \sqrt{2} |\xi|)}$
	& $\varrho \in (0,\infty)$
	& $d \geq 1$
	& $\gamma_0=\gamma_{\infty}=1$, $m=0$, $\ell=1$ \B\T \\ \hline
19
	& 
	& $\displaystyle \varrho |\xi| \coth((2 |\xi|)^{-1})- \varrho |\xi|^2$
	& $\varrho \in (0,\infty)$
	& $d \geq 1$
	& $\gamma_0=1$, $\gamma_{\infty}=2$, $m=0$, $\ell=1$ \B\T \\ \hline
20
	& 
	& $\displaystyle\varrho \log(\sinh(\sqrt{2} |\xi|))- \varrho \log(\sqrt{2} |\xi|)$
	& $\varrho \in (0,\infty)$
	& $d \geq 1$
	& $\gamma_0=2$, $\gamma_{\infty}=1$, $m=0$, $\ell=1$ \B\T \\ \hline
21
	& isotropic Meixner
	& $\displaystyle \varrho \log(\cosh(\sqrt{2} |\xi|))$
	& $\varrho \in (0,\infty)$
	& $d \geq 1$
	& $\gamma_0=2$, $\gamma_{\infty}=1$, $m=0$, $\ell=1$ \B\T \\ \hline
22
	& 
	& $\displaystyle |\xi| \log(1+\varrho \tanh(b |\xi|))$
	& $b \in (0,\infty)$, $\varrho \in (0,\infty)$
	& $d \geq 1$
	& $\gamma_0=\gamma_{\infty}=1$, $m=0$, $\ell=1$ \B\T \\ \hline
23
	& 
	& $\displaystyle \frac{\Gamma(\varrho |\xi|^2+1/2)}{\Gamma(\varrho |\xi|^2)}$
	& $\varrho \in (0,\infty)$
	& $d \geq 1$
	& $\gamma_0=2$, $\gamma_{\infty}=1$, $m \in (0,1/\sqrt{2\varrho})$, $\ell=1$ \B\T \\ \hline
24
	& 
	& $\displaystyle |\xi|^2 \frac{\Gamma(\alpha |\xi|^2+1-\alpha)}{\Gamma(\alpha |\xi|^2+1)}$
	& $\alpha \in (0,1)$
	& $d \geq 1$
	& $\gamma_0=2$, $\gamma_{\infty}=2-\alpha$, $\ell=1$, $m \in (0,\sqrt{\alpha^{-1}-1})$ \B\T \\ \hline
25
	& 
	& $\displaystyle\frac{\Gamma(\alpha |\xi|^2+1)}{\Gamma(\alpha |\xi|^2+1-\alpha)} - \frac{1}{\Gamma(1-\alpha)}$
	& $\alpha \in (0,1)$
	& $d \geq 1$
	& $\gamma_0=2$, $\gamma_{\infty}=\alpha$, $m \in (0,\sqrt{\alpha^{-1}})$ , $\ell=1$ \B\T \\ \hline
26
	&  Lamperti stable
	&  $\displaystyle \frac{\Gamma(|\xi|^2+\alpha+\varrho)}{\Gamma(|\xi|^2+\varrho)} - \frac{\Gamma(\alpha+\varrho)}{\Gamma(\varrho)}$ 
	& $\alpha \in (0,1)$, $\varrho \in (0,\infty)$
	& $d \geq 1$
	& $\gamma_0=2$, $\gamma_{\infty}=2\alpha$, $m \in (0,\sqrt{\varrho+\alpha})$, $\ell=1$ 
\end{tabular} \caption{(cont.)}
\end{table} 
\end{landscape}

If $(L_t)_{t \geq 0}$ is a subordinator, i.\,e.\ a L\'evy process with non-decreasing sample paths, then we can relax the assumptions of Theorem~\ref{levy-3}; the reason is that we know that the support of $L_t$ is contained in $[0,\infty)$ and therefore we have to establish upper bounds for $p_t(x)$ only for $x \geq 0$. Note that the characteristic exponent $\psi$ of a subordinator with Laplace exponent $f$ is given by $\psi(\xi) = f(-i \xi)$, $\xi \in \mbb{R}$.

\begin{kor} \label{levy-9}
	Let $(S_t)_{t \geq 0}$ be a subordinator with Laplace exponent $f$ satisfying \eqref{S2}, \eqref{S3}.
\begin{enumerate}[label*=\upshape (S\arabic*),ref=\upshape S\arabic*] 
			\item\label{S2} There exist $\theta \in (0,\tfrac{\pi}{2})$ and $m \geq 0$ such that $f$ has a holomorphic extension $F$ to \begin{equation*}
				\Upsilon := \Upsilon(m,\theta) := \{z \in \mbb{C}; -m<\re z<0\} \cup \{z \in \mbb{C} \backslash \{0\}; \arg z \in (\pi/2,\pi-\theta) \cup (-\pi+\theta,-\pi/2)\}.
			\end{equation*}
		\item\label{S3} There exist constants $c_1,c_2>0$ and $\gamma_0,\gamma_{\infty} \in (0,2]$ and a slowly varying (at infinity) increasing function $\ell: (0,\infty) \to (0,\infty)$ such that \begin{equation*}
				\re F(z) \geq ´\frac{c_1}{\ell(|z|)} |z|^{\gamma_{\infty}} \fa z \in \Upsilon, \; |z| \geq 1 
			\end{equation*}
		and \begin{equation*}
				|F(z)| \leq c_2 \ell(|z|) \left(|z|^{\gamma_0} \I_{\{|z| \leq 1\}} + |z|^{\gamma_{\infty}} \I_{\{|z| >1\}}\right) \fa z \in \Upsilon. 
			\end{equation*}
	\end{enumerate}
	Then $S_t$ has a density $p_t$ with respect to Lebesgue measure, \begin{equation*}
			p_t(x) := \frac{1}{2\pi} \int_{\mbb{R}} e^{-ix \cdot \xi} e^{-t f(-i\xi)} \, d\xi.
	\end{equation*}
	The density is infinitely often differentiable and satisfies the estimates \begin{align*}
			|p_t(x)| &\leq C \I_{[0,\infty)}(x) S(x,t) \\
			\left| \frac{\partial}{\partial t} p_t(x) \right| &\leq C \I_{[0,\infty)}(x) t^{-1} S(x,t)
		\end{align*}
	for any $x \in \mbb{R}$, $t \in (0,T]$ where $C=C(T)>0$ is an absolute constant and \begin{equation*}
			S(x,t) := S_m(x,t) := \exp \left( - \frac{m}{4} |x| \right) (1+\ell(ct^{-1/\gamma_{\infty}})) \begin{cases} t^{-1/\gamma_{\infty}}, & |x| \leq t^{1/\gamma_{\infty}} \wedge 1, \\  t/|x|^{1+\gamma_{\infty}}, & t^{1/\gamma_{\infty}} < |x| \leq 1, \\ t/|x|^{1+\gamma_{\infty} \wedge \gamma_0}, & |x|>1\end{cases}
		\end{equation*}
		for some absolute constant $c=c(T)>0$. 	Moreover, for any $k \in \mbb{N}$ and $T>0$ there exists a constant $c>0$ such that \begin{equation*}
			\left| \frac{\partial^k}{ \partial x^k} p_t(x) \right| \leq c \I_{[0,\infty)}(x) t^{-k/\gamma_{\infty}} S(x,t) \fa x \in \mbb{R}, t \in (0,T]. 
	\end{equation*}
\end{kor}

Corollary~\ref{levy-9} follows from the proof of Theorem~\ref{levy-3}, see \cite[pp.\ 131]{diss} or \cite[Section 4.9]{matters}. In Table~\ref{tab-sub} we have collected examples of Laplace exponents satisfying the assumptions of Corollary~\ref{levy-9}.

\begin{landscape} 
\begin{table}
\centering
	\begin{tabular}{p{20pt}p{230pt}p{110pt}p{190pt}}
		 & \textbf{Laplace Exponent} & \textbf{Parameter} & \textbf{Heat Kernel Estimate (Cor.~\ref{levy-9})}  \\
	\hline
		1
			& $\displaystyle \lambda^{\alpha}$  
			& $\alpha \in (0,1]$
			&$\gamma_0=\gamma_{\infty}=\alpha$, $m=0$, $\ell=1$ \B\T \\
		\hline
		2
			& $\displaystyle(\lambda+\varrho)^{\alpha}-\varrho^{\alpha}$
			& $\alpha \in (0,1)$, $\varrho \in (0,\infty)$
			&$\gamma_0=2$, $\gamma_{\infty}=\alpha$, $m \in (0,\varrho)$, $\ell=1$ \B\T \\ \hline
3
	& $\displaystyle\frac{\lambda}{\sqrt{\lambda+\varrho}}$
	& $\varrho \in (0,\infty)$
	& $\gamma_0=1$, $\gamma_{\infty}=1/2$, $m \in (0,\varrho)$, $\ell=1$ \B\T \\ \hline
4
	& $\displaystyle\frac{\lambda}{(\lambda+\varrho)^{\alpha}}$
	& $\alpha \in (0,1)$, $\varrho \in (0,\infty)$
	& $\gamma_0=1$, $\gamma_{\infty}=1-\alpha$, $m \in (0,\sqrt{\varrho})$, $\ell=1$\B\T \\ \hline
5
	& $\displaystyle\frac{\lambda^{\beta}-1}{\lambda^{\alpha}-1}-1$ (extended by continuity at $\lambda=1$)
	& $0<\alpha<\beta < 1$
	& $\gamma_0=\alpha$, $\gamma_{\infty}=\beta-\alpha$, $m=0$, $\ell=1$\B\T \\ \hline
6
	& $\displaystyle-\frac{\lambda^{\alpha}-1}{\lambda^{\alpha-1}-1}$ (extended by continuity at $\lambda=1$)
	& $\alpha \in (0,1)$
	& $\gamma_0=1-\alpha$, $\gamma_{\infty}=1$, $m =0$, $\ell=1$ \B\T \\ \hline
7
	& $\displaystyle\frac{\lambda^{\alpha}-1}{\lambda^{\alpha-1}-1}-1$ (extended by continuity at $\lambda=1$)
	& $\alpha \in (1,2]$
	& $\gamma_0=\alpha-1$, $\gamma_{\infty}=1$, $m=0$, $\ell=1$ \B\T \\ \hline
8
	& $\displaystyle \lambda \frac{\lambda^{\alpha}-\varrho^{\alpha}}{\lambda-\varrho}$ (extended by continuity at $\lambda=\varrho$)
	&$\alpha \in (0,1)$, $\varrho \in (0,\infty)$
	& $\gamma_0=1$, $\gamma_{\infty}=\alpha$, $m=0$, $\ell=1$\B\T \\ \hline
9
	& $\displaystyle(\lambda^{-\alpha}+\lambda^{-\beta})^{-1}$ (extended by continuity at $\lambda=0$)
	& $\alpha,\beta \in (0,1]$
	& $\gamma_0=\alpha \vee \beta$, $\gamma_{\infty}=\alpha \wedge \beta$, $m=0$, $\ell=1$ \B\T \\ \hline
10
	& $\displaystyle \sqrt{\lambda} (1-e^{-2 \varrho \sqrt{\lambda}})$
	& $\varrho \in (0,\infty)$
	& $\gamma_0=1$, $\gamma_{\infty}=1/2$, $m=0$, $\ell=1$ \B\T \\ \hline
11
	& $\displaystyle \sqrt{\lambda} (1+e^{-2 \varrho \sqrt{\lambda}})$
	& $\varrho \in (0,\infty)$
	& $\gamma_0=\gamma_{\infty}=1/2$, $m=0$, $\ell=1$ \B\T \\ \hline
12
	& $\displaystyle\varrho \lambda (\lambda+1) \log(1+1/\lambda)$ %|\xi|^2 \log(1+|\xi|^{-2}) \to 1$ as $|\xi| \to \infty
	& $\varrho \in (0,\infty)$
	& $\gamma_0=\gamma_{\infty}=1$, $m=0$, $\ell(r)=\log(r \vee e)$ \end{tabular} \caption{Examples of Laplace exponents satisfying the assumptions of Corollary~\ref{levy-9}} \label{tab-sub}. \end{table}
	\setcounter{table}{1}
\begin{table} \centering
\begin{tabular}{p{20pt}p{230pt}p{110pt}p{190pt}}
	& \textbf{Laplace Exponent} & \textbf{Parameter} & \textbf{Heat Kernel Estimate (Cor.~\ref{levy-9})}  \\ \hline
13
	& $\displaystyle \varrho \frac{\lambda (\lambda+1)}{(\lambda+2) \log(\lambda+2)}$
	& $\varrho \in (0,\infty)$
	& $\gamma_0=\gamma_{\infty}=1$, $m \in (0,2)$, $\ell(r)=\log(r \vee e)$ \B\T \\ \hline
14
	& $\displaystyle \sqrt{\lambda} \arctan(\varrho \sqrt{\lambda})$
	& $\varrho \in (0,\infty)$
	& $\gamma_0=\gamma_{\infty}=1/2$, $m=0$, $\ell=1$ \B\T \\ \hline
15
	& $\displaystyle (\lambda+\varrho)^{\alpha} \cos \big(\alpha \arctan \sqrt{\tfrac{\lambda}{\varrho}}\big)-\varrho^{\alpha}$
	& $\alpha \in (0,1)$, $\varrho \in (0,\infty)$
	& $\gamma_0=1$, $\gamma_{\infty}=\alpha$, $m \in (0,\varrho)$, $\ell=1$ \B\T \\ \hline
16
	& $\displaystyle \varrho \sqrt{\lambda} \frac{\cosh^2(\sqrt{2\lambda})}{\sinh(2 \sqrt{2\lambda})}$
	& $\varrho \in (0,\infty)$
	& $\gamma_0=\gamma_{\infty}=1/2$, $m=0$, $\ell=1$ \B\T \\ \hline
17
	& $\displaystyle \varrho \sqrt{\lambda} \frac{\sinh^2(\sqrt{2\lambda})}{\sinh(2 \sqrt{2\lambda})}$
	& $\varrho \in (0,\infty)$
	& $\gamma_0=\gamma_{\infty}=1/2$, $m=0$, $\ell=1$ \B\T \\ \hline
18
	& $\displaystyle \varrho \sqrt{\lambda} \coth((2\sqrt{\lambda})^{-1})- \varrho \lambda$
	& $\varrho \in (0,\infty)$
	& $\gamma_0=1/2$, $\gamma_{\infty}=1$, $m=0$, $\ell=1$ \B\T \\ \hline
19
	& $\displaystyle\varrho \log(\sinh(\sqrt{2\lambda}))- \varrho \log(\sqrt{2\lambda})$
	& $\varrho \in (0,\infty)$
	& $\gamma_0=1$, $\gamma_{\infty}=1/2$, $m=0$, $\ell=1$ \B\T \\ \hline
20
	& $\displaystyle \varrho \log(\cosh(\sqrt{2\lambda}))$
	& $\varrho \in (0,\infty)$
	& $\gamma_0=1$, $\gamma_{\infty}=1/2$, $m=0$, $\ell=1$ \B\T \\ \hline
21
	& $\displaystyle \sqrt{\lambda} \log(1+\varrho \tanh(b \sqrt{\lambda}))$
	& $b \in (0,\infty)$, $\varrho \in (0,\infty)$
	& $\gamma_0=\gamma_{\infty}=1/2$, $m=0$, $\ell=1$ \B\T \\ \hline
22
	& $\displaystyle \frac{\Gamma(\varrho \lambda+1/2)}{\Gamma(\varrho \lambda)}$
	& $\varrho \in (0,\infty)$
	& $\gamma_0=1$, $\gamma_{\infty}=1/2$, $m \in (0,1/(2\varrho))$, $\ell=1$ \B\T \\ \hline
23
	& $\displaystyle \lambda \frac{\Gamma(\alpha \lambda+1-\alpha)}{\Gamma(\alpha \lambda+1)}$
	& $\alpha \in (0,1)$
	& $\gamma_0=1$, $\gamma_{\infty}=1-\alpha$, $\ell=1$, $m \in (0,\alpha^{-1}-1)$ \B\T \\ \hline
24
	& $\displaystyle\frac{\Gamma(\alpha \lambda+1)}{\Gamma(\alpha \lambda+1-\alpha)} - \frac{1}{\Gamma(1-\alpha)}$
	& $\alpha \in (0,1)$
	& $\gamma_0=1$, $\gamma_{\infty}=\alpha$, $m \in (0,\alpha^{-1})$ , $\ell=1$ \B\T \\ \hline
25
	&  $\displaystyle \frac{\Gamma(\lambda+\alpha+\varrho)}{\Gamma(\lambda+\varrho)} - \frac{\Gamma(\alpha+\varrho)}{\Gamma(\varrho)}$ 
	& $\alpha \in (0,1)$, $\varrho \in (0,\infty)$
	& $\gamma_0=1$, $\gamma_{\infty}=\alpha$, $m \in (0,\varrho+\alpha)$, $\ell=1$ 
\end{tabular} \caption{(cont.)}
\end{table} 
\end{landscape}

\section{Existence result for L\'evy-type processes} \label{ltp}

Let $q: \mbb{R}^d \times \mbb{R}^d \to \mbb{C}$ be defined by \begin{equation*}
	q(x,\xi) = \psi_{\alpha(x)}(\xi), \qquad x,\xi \in \mbb{R}^d
\end{equation*}
for a family $(\psi_{\beta})_{\beta \in I}$ of continuous negative definite functions, a H\"{o}lder continuous mapping $\alpha: \mbb{R}^d \to I$ and a set of parameters $I \subseteq \mbb{R}^n$. Our main result, Theorem~\ref{ltp-3}, gives a sufficient condition on $(\psi_{\beta})_{\beta \in I}$ for the existence of a
rich L\'evy-type process with symbol $q$. 

\begin{thm} \label{ltp-3}
	Let $I \subseteq \mbb{R}^n$ be open and convex and $m \geq 0$.  Let $(\psi_{\beta})_{\beta \in I}$ be a family of continuous negative definite functions $\psi_{\beta}:\mbb{R}^d \to \mbb{C}$ with $\psi_{\beta}(0)=0$ for all $\beta \in I$. Suppose that there exist $\theta \in (0,\pi/2)$ and constants $c_1,c_2,c_3>0$ such that each $\psi_{\beta}$, $\beta \in I$, satisfies \eqref{LTP1}--\eqref{LTP4}.
	\begin{description}
		\item[\normalfont{(LTP1)}]\label{LTP1} $\psi_{\beta}$ is rotationally invariant for each $\beta \in I$, i.\,e.\ there exists $\Psi_{\beta}: \mbb{R} \to \mbb{R}$ such that $\psi_{\beta}(\xi) = \Psi_{\beta}(|\xi|)$, $\xi \in \mbb{R}^d$. If $m>0$: $\Psi_{\beta}(r)=\Psi_{\beta}(-r)$ for all $r \geq 0$.
		\item[\normalfont{(LTP2)}]\label{LTP2} $\psi_{\beta}$ has a holomorphic extension to the domain $\Omega=\Omega(m,\theta)$ defined in \eqref{domain}.
		\item[\normalfont{(LTP3)}]\label{LTP3} There exist a measurable mapping $\gamma_0: I \to (0,2]$, a H\"{o}lder continuous mapping $\gamma_{\infty}: I \to (0,2]$ and a slowly varying (at infinity) increasing function $\ell: (0,\infty) \to (0,\infty)$ such that \begin{equation*}
			\re \Psi_{\beta}(z) \geq \frac{c_1}{\ell(|z|)} |\re z|^{\gamma_{\infty}(\beta)} \fa z \in \Omega, |z| \geq 1, \beta \in I
		\end{equation*}
		and \begin{equation*}
			|\Psi_{\beta}(z)| \leq c_2 \ell(|z|) \left(|z|^{\gamma_0(\beta)} \I_{\{|z| \leq 1\}} + |z|^{\gamma_{\infty}(\beta)} \I_{\{|z|>1\}}\right) \fa z \in \Omega, \; \beta \in I.
		\end{equation*}
		Moreover, $\gamma_{\infty}^L := \inf_{\beta \in I} \gamma_{\infty}(\beta)>0$, $\gamma_0^L := \inf_{\beta \in I} \gamma_0(\beta)>0$.
		\item[\normalfont{(LTP4)}]\label{LTP4} The partial derivative $\frac{\partial}{\partial \beta_j} \Psi_{\beta}(r)$ exists for all $r \in \mbb{R}$ and extends holomorphically to $\Omega$ for all $j \in \{1,\ldots,n\}$ and $\beta \in I$. Moreover, \begin{align*}
			\left| \frac{\partial}{\partial \beta_j} \Psi_{\beta}(z) \right| &\leq c_3 (1+\ell(|z|)) \left(|z|^{\gamma_0(\beta)} \I_{\{|z| \leq 1\}} + |z|^{\gamma_{\infty}(\beta)} \I_{\{|z|>1\}}\right) \fa z \in \Omega, \; \beta \in I.
		\end{align*}
		\end{description}
	Then for any H\"{o}lder continuous mapping $\alpha: \mbb{R}^d \to I$ there exists a rich L\'evy-type process $(X_t)_{t \geq 0}$ with symbol \begin{equation*}
		q(x,\xi) := \psi_{\alpha(x)}(\xi), \qquad x,\xi \in \mbb{R}^d.
	\end{equation*}
\end{thm}

To prove Theorem~\ref{ltp-3} we use the parametrix method, cf.\ \cite[Chapter 4]{diss} and \cite[Chapter 4]{matters}. The idea is to construct the transition density as the fundamental solution of the Cauchy problem for the operator $(\partial_t-L)$ where $L$ equals, when restricted to $C_c^{\infty}(\mbb{R}^d)$, the pseudo-differential operator \begin{equation*}
	Lf(x) = - \int_{\mbb{R}^d} q(x,\xi) e^{ix \cdot \xi} \hat{f}(\xi) \, d\xi, \qquad x \in \mbb{R}^d.
\end{equation*}
The parametrix method gives a candidate for the fundamental solution, and the main part of the proof is to verify that this candidate is indeed a fundamental solution to the Cauchy problem and the transition density of a Feller process. \par
As a by-product of the parametrix construction, we get the following additional information on $(X_t)_{t \geq 0}$ and its transition semigroup $(P_t)_{t \geq 0}$.

\begin{thm} \label{ltp-5}
	Under the assumptions of Theorem~\ref{ltp-3}, the L\'evy-type process $(X_t)_{t \geq 0}$ with symbol $q(x,\xi) = \psi_{\alpha(x)}(\xi)$ has the following additional properties:
 	\begin{enumerate}
		\item\label{ltp-5-i} The associated semigroup $(P_t)_{t \geq 0}$ has the strong Feller property, i.\,e.\ $P_t f \in C_b(\mbb{R}^d)$ for any $f \in \mc{B}_b(\mbb{R}^d)$.
		\item\label{ltp-5-ii} $C_c^{\infty}(\mbb{R}^d)$ is a core for the generator $L$ and $C_{\infty}^2(\mbb{R}^d) \subseteq \mc{D}(L)$, \begin{equation*}
			Lf(x) = b(x) \cdot \nabla f(x) + \frac{1}{2} \tr(Q(x) \cdot \nabla^2 f(x)) + \int_{y \neq 0}(f(x+y)-f(x)-\nabla f(x) \cdot y \I_{B(0,1)}(y)) \, \nu(x,dy)
		\end{equation*}
		for any $f \in C_{\infty}^2(\mbb{R}^d)$; here $(b(x),Q(x),\nu(x,dy))$ denotes the L\'evy triplet associated with $q(x,\cdot)$. There exists a constant $C>0$ such that \begin{equation*}
			\|Lf\|_{\infty} \leq C \sum_{0 \leq |\alpha| \leq 2} \|\partial^{\alpha} f\|_{\infty} =  C \|f\|_{(2)} \fa f \in C_{\infty}^2(\mbb{R}^d).
		\end{equation*}
		Moreover, $P_t f \in \mc{D}(L)$ for all $t>0$ and $f \in C_{\infty}(\mbb{R}^d)$.
		\item\label{ltp-5-iii} The distribution $\mbb{P}^x(X_t \in \cdot)$ has a density $p(t,x,\cdot)$ with respect to Lebesgue measure for all $t>0$ and $x \in \mbb{R}^d$. The mapping $p:(0,\infty) \times \mbb{R}^d \times \mbb{R}^d \to [0,\infty)$ is continuous and differentiable with respect to $t$.
		\item\label{ltp-5-iv} The transition density $p$ is a fundamental solution to the Cauchy problem for the operator $(\partial_t-L)$, i.\,e.\ $p(t,\cdot,y)$ converges weakly to $\delta_x$ as $t \to 0$, $(0,\infty) \ni t \mapsto p(t,x,y)$ is differentiable, $p(t,\cdot,y) \in \mc{D}(L)$ for all $t>0$, $y \in \mbb{R}^d$ and \begin{equation*}
			(\partial_t-L_x) p(t,x,y)=0 \fa t>0, x,y \in \mbb{R}^d.
		\end{equation*}
		\item\label{ltp-5-v} The $(L,C_c^{\infty}(\mbb{R}^d))$-martingale problem is well-posed; its unique (in the sense of finite-dimensional distributions) solution is $(X_t)_{t \geq 0}$.
		\item\label{ltp-5-vi} Denote by $\varrho \in (0,1]$ the H\"{o}lder exponent of $\alpha$ and choose $\gamma \in (0,1/\gamma_{\infty}^U]$ such that $\kappa := \gamma \min\{\varrho, (-d+\gamma_{\infty}^U)+1\}>0$. Define \begin{equation*}
			S(x,\beta,t) := \exp \left( - \frac{m}{4} |x| \right) \begin{cases} t^{-d/\gamma_{\infty}(\beta)}, & |x| \leq t^{1/\gamma_{\infty}(\beta)} \wedge 1, \\ t/|x|^{d+\gamma_{\infty}(\beta)}, & t^{1/\gamma_{\infty}(\beta)} \leq |x| \leq 1, \\ t/|x|^{d+\gamma_0(\beta) \wedge \gamma_{\infty}(\beta)}, &|x|>1. \end{cases}
		\end{equation*}
		For any $T>0$ there exists a constant $C=C(T)>0$ such that\begin{align*}
			|p(t,x,y)| \leq C S(x-y,\alpha(y),t) + Ct^{\kappa} \frac{1}{1+|x-y|^{\gamma_0^L \wedge \gamma_{\infty}^L}} \exp \left( - \frac{m}{4} |x-y| \right) \\ 
			|\partial_t p(t,x,y)| \leq C t^{-1} S(x-y,\alpha(y),t) + Ct^{-1+\kappa} \frac{1}{1+|x-y|^{\gamma_0^L \wedge \gamma_{\infty}^L}} \exp \left( - \frac{m}{4} |x-y| \right)
		\end{align*} 
		for all $x,y \in \mbb{R}^d$	and $t \in (0,T]$. 
	\end{enumerate}
\end{thm}

\begin{bems} \label{ltp-6} \begin{enumerate}
	\item In dimension $d=1$ we can drop the assumption \eqref{LTP1} of rotational invariance, see Theorem~\ref{ltp-7} below.
	\item The constant $m/4$ in the definition of $S$, cf.\ Theorem~\ref{ltp-5}\ref{ltp-5-vi}, may be replaced by $m(1-\delta)$ for any $\delta \in (0,1)$; see also Remark~\ref{levy-5}.
	\item In Theorem~\ref{ltp-3} we make separate assumptions on the regularity of $I \ni \beta \mapsto \psi_{\beta}(\xi)$ (differentiability) and $\mbb{R}^d \ni x \mapsto \alpha(x) \in I$ (H\"{o}lder continuity). Note that this is much weaker than assuming differentiability of $x \mapsto q(x,\xi) = \psi_{\alpha(x)}(\xi)$. For instance, if $\psi_{\beta}(\xi) := |\xi|^{\beta}$, then the assumptions of Theorem~\ref{ltp-3} are satisfied for any H\"{o}lder continuous function $\alpha$; in contrast, differentiability of $x \mapsto q(x,\xi) = |\xi|^{\alpha(x)}$ requires differentiability of $\alpha$.
	\item The first order approximation of the transition probability $p$ is given by \begin{equation*}
		p_0(t,x,y) := \frac{1}{(2\pi)^d} \int_{\mbb{R}^d} e^{-i(x-y) \cdot \xi} e^{-t \psi_{\alpha(y)}(\xi)} \, d\xi; 
	\end{equation*}
	it is possible to derive upper bounds for $|p(t,x,y)-p_0(t,x,y)|$, cf.\ \cite[Theorem 3.8]{matters} for details.
\end{enumerate} \end{bems}

\begin{thm}[Case $d=1$] \label{ltp-7}
	Let $I \subseteq \mbb{R}^n$ be an open convex set and $m \geq 0$. Suppose that $(\psi_{\beta})_{\beta \in I}$ is a family of continuous negative definite functions $\psi_{\beta}:\mbb{R} \to \mbb{C}$, $\psi_{\beta}(0)=0$, such that $\Psi_{\beta}(\xi) := \psi_{\beta}(\xi)$, $\beta \in I$, satisfies \eqref{LTP2}-\eqref{LTP4}. Then the statements of Theorem~\ref{ltp-3} and Theorem~\ref{ltp-5} remain valid; in particular, there exists a rich L\'evy-type process $(X_t)_{t \geq 0}$ with symbol $q(x,\xi) =\psi_{\alpha(x)}(\xi)$, $x,\xi \in \mbb{R}$, for any H\"{o}lder continuous mapping $\alpha: \mbb{R} \to I$.
\end{thm}

The next theorem shows that in dimension $d=1$ the transition probability $p(t,x,y)$ is differentiable with respect to $x$ provided that the mappings $I \ni \beta \mapsto \psi_{\beta}(\xi)$ and $\alpha: \mbb{R}^d \to I$ are sufficiently smooth.

\begin{thm} \label{ltp-9}
	Let $(\psi_{\beta})_{\beta \in I}$ be as in Theorem~\ref{ltp-7} and assume additionally that there exists a constant $c_4>0$ such that \eqref{LTP5} holds. 
	\begin{description}  
		\item[\normalfont{(LTP5)}]\label{LTP5} $\frac{\partial^2}{\partial \beta_j^2} \psi_{\beta}(\xi)$ exists for all $\xi \in \mbb{R}$, $j \in \{1,\ldots,n\}$ and has a holomorphic extension to $\Omega$ satisfying \begin{align*}
			\left| \frac{\partial^2}{\partial \beta_j^2} \psi_{\beta}(z) \right|
			&\leq c_4 (1+\ell(|z|))( |z|^{\gamma_{0}(\beta)} \I_{\{|z| \leq 1\}} + |z|^{\gamma_{\infty}(\beta)} \I_{\{|z|>1\}}) , \qquad z \in \Omega, \beta \in I
		\end{align*}
		where $\ell$ denotes the slowly varying function from \eqref{LTP3}.
	\end{description}
	Let $\alpha: \mbb{R} \to I$ be such that $\alpha \in C_b^2(\mbb{R})$. Denote by $(X_t)_{t \geq 0}$ the L\'evy-type process from Theorem~\ref{ltp-7} with symbol $q(x,\xi) = \psi_{\alpha(x)}(\xi)$ and transition density $p$. Then: \begin{enumerate}
		\item The transition probability $p(t,x,y)$ is continuously differentiable with respect to $x$ for any $t>0$ and $y \in \mbb{R}$. For any $T>0$ there exists a constant $C=C(T)>0$ such that \begin{equation*}
			\left| \frac{\partial}{\partial x} p(t,x,y) \right| \leq Ct^{-1/\gamma_{\infty}^L} \left[S(x-y,\alpha(y),t) + t^{\kappa} \frac{1}{1+|x-y|^{d+\gamma_0^L \wedge \gamma_{\infty}^L}} \exp \left(-\frac{m}{4} |x-y| \right) \right]
		\end{equation*}
		for all $t \in (0,T]$ and $x,y \in \mbb{R}$; see Theorem~\ref{ltp-5}\,\ref{ltp-5-vi} for the definition of $\kappa$ and $S$.
		\item The semigroup $(P_t)_{t \geq 0}$ asssociated with the L\'evy-type process $(X_t)_{t \geq 0}$ satisfies the gradient estimate	\begin{equation*}
			\sup_{x \in \mbb{R}} \left| \frac{\partial}{\partial x} P_t f(x) \right| \leq C t^{-1/\gamma_{\infty}^L} \|f\|_{\infty} \fa  t \in (0,T], \; f \in \mc{B}_b(\mbb{R})
		\end{equation*}
		for some absolute constant $C=C(T)>0$.
		\item Suppose additionally that each $\psi_{\beta}: \mbb{R} \to \mbb{R}$, $\beta \in I$, is even. Then for any $T>0$ there exist constants $C_1,C_2,C_3>0$ such that \begin{equation*}
			p(t,x,y) \geq C_1 t^{-1/\gamma_{\infty}(\alpha(y))} \big( 1- C_2 t^{-1/\gamma_{\infty}(\alpha(y))} |x-y| - C_3 t^{\kappa} \big)^+ 
		\end{equation*}
		for all $x,y \in \mbb{R}$, $t \in (0,T]$.
		\item If $\psi_{\beta}: \mbb{R} \to \mbb{R}$ is an even function for all $\beta \in I$, then $(X_t)_{t \geq 0}$ is $\lambda$-irreducible, i.\,e. \begin{equation*}
			\int_{(0,\infty)} \mbb{P}^x(X_t \in B) \, dt > 0 
		\end{equation*}
		for all $x \in \mbb{R}$ and $B \in \mc{B}(\mbb{R})$ with $\lambda(B)>0$.
	\end{enumerate}
\end{thm}

The remaining part of this article is devoted to applications of the above results. First, we state an existence result for ``stable-like'' processes. The most popular examples are isotropic stable-like processes (processes with variable index of stability); this corresponds to symbols of the form $q(x,\xi) = |\xi|^{\alpha(x)}$. Bass \cite{bass} proved the well-posedness of the associated martingale problem in dimension $d=1$ for Dini continuous functions $\alpha$, and, more recently, Kolokoltsov \cite{kol} established the existence of Feller processes with symbol $q(x,\xi) = |\xi|^{\alpha(x)}$ in dimension $d \geq 1$ for H\"{o}lder continuous mappings $\alpha$. Using the results from the first part of this Section, we can derive existence results for many stable-like processes, for instance relativistic stable-like and Lamperti stable-like processes. With the exception of the well-studied isotropic stable case, such existence results were so far only known under much stronger regularity assumptions; e.\,g.\ a general existence result by Hoh \cite{hoh_hab} requires $\alpha \in C^{5d+3}(\mbb{R}^d)$.

\begin{bsp} \label{ltp-11}
	Let $q(x,\xi)$ be one of the following functions.
    \begin{enumerate}
	\item
        (isotropic stable-like) $q(x,\xi) = |\xi|^{\alpha(x)}$ where $\alpha: \mbb{R}^d \to (0,2]$ is a H\"{o}lder continuous mapping such that $\inf_{x \in \mbb{R}^d} \alpha(x)>0$.
    \item
        (relativistic stable-like) $q(x,\xi) = (|\xi|^2+\varrho(x)^2)^{\alpha(x)/2}-\varrho(x)^{\alpha(x)}$ for H\"{o}lder continuous mappings $\alpha: \mbb{R}^d \to (0,2)$ and $\varrho: \mbb{R}^d \to (0,\infty)$ such that
        \begin{equation*}
		    \inf_{x \in \mbb{R}^d} \alpha(x)>0
            \quad\text{and}\quad
            0 < \varrho^L := \inf_{x \in \mbb{R}^d} \varrho(x) \leq \sup_{x \in \mbb{R}^d} \varrho(x) < \infty.
	    \end{equation*}
    \item
        (Lamperti stable-like) $q(x,\xi) = (|\xi|^2+\varrho(x))_{\alpha(x)}-(\varrho(x))_{\alpha(x)}$ for H\"{o}lder continuous mappings $\alpha: \mbb{R}^d \to (0,1)$ and $\varrho: \mbb{R}^d \to (0,\infty)$ such that
        \begin{equation*}
            0<\inf_{x \in \mbb{R}^d} \alpha(x) \leq \sup_{x \in \mbb{R}^d} \alpha(x)<1
            \quad \text{and} \quad
            0 < \varrho^L := \inf_{x \in \mbb{R}^d} \varrho(x) \leq \sup_{x \in \mbb{R}^d} \varrho(x) < \infty;
		\end{equation*}
		here $(r)_{\alpha} := \Gamma(r+\alpha)/\Gamma(r)$ denotes the Pochhammer symbol.
	\item
        (TLP-like)
        $q(x,\xi) = (|\xi|^2+\varrho(x)^2)^{\alpha(x)/2} \cos\big[\alpha(x) \arctan \tfrac{|\xi|}{\varrho(x)}\big]-\varrho(x)^{\alpha(x)}$ for H\"{o}lder continuous mappings $\alpha: \mbb{R}^d \to (0,1)$ and $\varrho: \mbb{R}^d \to (0,\infty)$ such that
        \begin{equation*}
            0<\inf_{x \in \mbb{R}^d} \alpha(x) \leq \sup_{x \in \mbb{R}^d} \alpha(x)<1
            \quad\text{and}\quad
            0 < \varrho^L := \inf_{x \in \mbb{R}^d} \varrho(x) \leq \sup_{x \in \mbb{R}^d} \varrho(x) < \infty.
		\end{equation*}
	\end{enumerate}
	By Theorem~\ref{ltp-3}, there exists a rich L\'evy-type process $(X_t)_{t \geq 0}$ with symbol $q$. The process $(X_t)_{t \geq 0}$ has the properties listed in Theorem~\ref{ltp-5}; the heat kernel estimate \ref{ltp-5}.\ref{ltp-5-vi} holds with \begin{enumerate}
		\item $m=0$, $\gamma_{\infty}(\alpha(x)) = \gamma_0(\alpha(x))=\alpha(x)$
		\item $m=\varrho^L$, $\gamma_{\infty}(\alpha(x),\varrho(x))=\alpha(x)$, $\gamma_0(\alpha(x),\varrho(x))=2$,
		\item $m=\sqrt{\varrho^L}$, $\gamma_{\infty}(\alpha(x),\varrho(x))=2\alpha(x)$, $\gamma_0(\alpha(x),\varrho(x))=2$,
		\item $m=\varrho^L$, $\gamma_{\infty}(\alpha(x),\varrho(x)) = \alpha(x)$, $\gamma_0(\alpha(x),\varrho(x))=2$.
	\end{enumerate}
	Let us remark that is possible to obtain further information on the richness of the domain of the infinitesimal generator of $(X_t)_{t \geq 0}$, cf.\ \cite[Example 4.11]{hke}. 
\end{bsp}

Applying Theorem~\ref{ltp-7} we obtain in a similar fashion an existence result for one-dimensional rich L\'evy-type processes with symbol \begin{equation*}
	q(x,\xi) = (\kappa(x)^2+(\xi-ib(x))^2)^{\alpha(x)/2} - (\kappa(x)^2-b(x)^2)^{\alpha(x)/2}, \qquad x,\xi \in \mbb{R}
\end{equation*}
for H\"{o}lder continuous bounded mappings $b: \mbb{R} \to \mbb{R}$, $\alpha: \mbb{R} \to (0,2)$, $\kappa: \mbb{R} \to (0,\infty)$ such that \begin{equation*}
		\alpha^L := \inf_{x \in \mbb{R}} \alpha(x)>0, \qquad \kappa^L := \inf_{x \in \mbb{R}} \kappa(x)>0, \qquad \kappa^L-\|b\|_{\infty}>0;
\end{equation*}
we call such a L\'evy-type process an NTS-like process; NTS is short for normal tempered stable. More generally, it is possible to consider symbols of the form \begin{equation*}
	q(x,\xi) = f_{\alpha(x)}(|\xi|^2)
\end{equation*}
for a family of Bernstein functions $(f_{\beta})_{\beta \in I}$, this leads to, so-called, variable order subordination; see \cite[Section 5.1]{matters} for a general existence result. \par
Further examples of families of continuous negative definite functions satisfying \eqref{LTP1}-\eqref{LTP4} are listed in Table~\ref{tab-ltp}; in Table~\ref{tab-ltp} we use $C^{>0}(I)$ to denote the space of bounded H\"{o}lder  continuous functions $f: \mbb{R}^d \to I$ satisfying \begin{equation*}
	f^L(x) := \inf_{x \in \mbb{R}^d} f(x) \in I \quad \text{and} \quad f^U(x) := \sup_{x \in \mbb{R}^d} f(x) \in I.
\end{equation*}

\begin{landscape} \begin{table} \begin{tabular}{llllll}
  \textbf{Name} & \textbf{Symbol} & \textbf{Assumptions} & \textbf{Dim.} & \textbf{Parameters (Thm.~\ref{ltp-3})} \\
\hline
	 isotropic $\alpha$-stable-like
	& $\displaystyle|\xi|^{\alpha(x)}$  
	& $\alpha \in C^{>0}((0,2])$
	& $d \geq 1$
	& $\gamma_{\infty}(\alpha(x)) = \gamma_0(\alpha(x))=\alpha(x)$, $m=0$ \B\T \\ \hline
	
	 relativistic stable-like
	& $\displaystyle(|\xi|^2+\varrho^2(x))^{\alpha(x)/2}-\varrho(x)^{\alpha(x)}$
	& \vtop{ \hbox{\strut $\alpha \in C^{>0}((0,2])$} \hbox{\strut $\varrho \in C^{>0}((0,\infty))$}}
	& $d \geq 1$
	& \vtop{ \hbox{\strut $\gamma_0(\alpha(x),\varrho(x))=2$, $m \in (0,\varrho^L)$} \hbox{\strut $\gamma_{\infty}(\alpha(x),\varrho(x))=\alpha(x)$}}\B\T \\ \hline

	 NTS-like
	& $\displaystyle(\kappa(x)^2+(\xi-ib(x))^2)^{\alpha(x)/2}-(\kappa(x)^2-b(x)^2)^{\alpha(x)/2}$
	& \vtop{ \hbox{\strut $\alpha \in C^{>0}((0,2])$} \hbox{\strut $\kappa \in C^{>0}((0,\infty))$} \hbox{\strut $b \in C^{>0}(\mbb{R})$} \hbox{\strut $\kappa^L-\|b\|_{\infty}>0$}}
	& $d=1$
	& \vtop{ \hbox{\strut $\gamma_0(\alpha(x),b(x),\kappa(x))=2$} \hbox{\strut $\gamma_{\infty}(\alpha(x),b(x),\kappa(x))=\alpha(x)$} \hbox{\strut $m \in (0,\kappa^L-\|b\|_{\infty})$}} \B\T \\ \hline

	& $\displaystyle\frac{|\xi|^2}{\sqrt{|\xi|^2+\varrho(x)}}$
	& $\varrho \in C^{>0}((0,\infty))$
	& $d \geq 1$
	& \vtop{\hbox{\strut $\gamma_0(\varrho(x))=2$, $\gamma_{\infty}(\varrho(x))=1$} \hbox{\strut $m \in (0,\sqrt{\varrho^L)}$}}\B\T \\ \hline

	& $\displaystyle\frac{|\xi|^2}{(|\xi|^2+\varrho(x))^{\alpha(x)}}$
	& \vtop{\hbox{\strut $\alpha \in C^{>0}((0,2))$} \hbox{\strut $\varrho \in C^{>0}((0,\infty))$}}
	& $d \geq 1$
	& \vtop{\hbox{\strut $\gamma_0(\alpha(x),\varrho(x))=2$, $m \in (0,\sqrt{\varrho^L})$} \hbox{\strut $\gamma_{\infty}(\alpha(x),\varrho(x))=2-\alpha(x)$}} \B\T \\ \hline

	& $\displaystyle\frac{|\xi|^{\beta(x)}-1}{|\xi|^{\alpha(x)}-1}-1$ (extended by continuity at $\xi=1$)
	& \vtop{\hbox{\strut $\alpha,\beta \in C^{>0}((0,1))$} \hbox{\strut $(\beta-\alpha)^L>0$}}
	& $d \geq 1$
	& \vtop{\hbox{\strut $\gamma_0(\alpha(x),\beta(x))=\alpha(x)$, $m=0$} \hbox{\strut $\gamma_{\infty}(\varrho(x))=\beta(x)-\alpha(x)$}} \B\T \\ \hline

	& $\displaystyle-\frac{|\xi|^{\alpha(x)}-1}{|\xi|^{\alpha(x)-2}-1}$ (extended by continuity at $\xi=1$)
	& $\alpha \in C^{>0}((0,2))$
	& $d \geq 1$
	& \vtop{\hbox{\strut $\gamma_0(\alpha(x))=2-\alpha(x)$, $\gamma_{\infty}(\alpha(x))=2$} \hbox{\strut $m =0$}}  \B\T \\ \hline

	& $\displaystyle\frac{|\xi|^{\alpha(x)}-1}{|\xi|^{\alpha(x)-2}-1}-1$ (extended by continuity at $\xi=1$)
	& $\alpha \in C^{>0}((2,4])$
	& $d \geq 1$
	& \vtop{\hbox{\strut $\gamma_0(\alpha(x))=\alpha(x)-2$, $\gamma_{\infty}(\alpha(x))=2$} \hbox{\strut $m=0$}}  \B\T \\ \hline

	& $\displaystyle|\xi|^2 \frac{|\xi|^{\alpha(x)}-\varrho(x)^{\alpha(x)}}{|\xi|^2-\varrho(x)^2}$ (extended by continuity at $\xi=\varrho(x)$)
	& \vtop{\hbox{\strut $\alpha \in C^{>0}((0,2))$} \hbox{\strut $\varrho \in C^{>0}((0,\infty))$}}
	& $d \geq 1$
	& \vtop{\hbox{\strut $\gamma_0(\alpha(x),\varrho(x))=2$, $m=0$} \hbox{\strut $\gamma_{\infty}(\alpha(x),\varrho(x))=\alpha(x)$}}\B\T \\ \hline

	& $\displaystyle(|\xi|^{-\alpha(x)}+|\xi|^{-\beta(x)})^{-1}$ (extended by continuity at $\xi=0$)
	& $\alpha,\beta \in C^{>0}((0,2])$
	& $d \geq 1$
	& \vtop{\hbox{\strut $\gamma_0(\alpha(x),\beta(x))=\alpha(x) \vee \beta(x)$, $m=0$} \hbox{\strut $\gamma_{\infty}(\alpha(x),\beta(x)) = \alpha(x) \wedge \beta(x)$}}  \B\T \\ \hline

	& $\displaystyle|\xi| (1-e^{-2 \varrho(x) |\xi|})$
	& $\varrho \in C^{>0}((0,\infty))$
	& $d \geq 1$
	& $\gamma_0(\varrho(x))=2$, $\gamma_{\infty}(\varrho(x))=1$, $m=0$ \B\T \\ \hline

	& $\displaystyle|\xi| (1+e^{-2 \varrho(x) |\xi|})$
	& $\varrho \in C^{>0}((0,\infty))$
	& $d \geq 1$
	& $\gamma_0(\varrho(x))=\gamma_{\infty}(\varrho(x))=1$, $m=0$ \B\T \\ \hline
	 
	& $\displaystyle\varrho(x) |\xi|^2 (|\xi|^2+1) \log(1+|\xi|^{-2})$ %|\xi|^2 \log(1+|\xi|^{-2}) \to 1$ as $|\xi| \to \infty
	& $\varrho \in C^{>0}((0,\infty))$
	& $d \geq 1$
	& $\gamma_0(\varrho(x))=\gamma_{\infty}(\varrho(x))=2$, $m=0$
\end{tabular} \caption{Examples of admissible symbols} \label{tab-ltp}  \end{table}
\setcounter{table}{2}
\begin{table} \begin{tabular}{lllll}
  \textbf{Name} & \textbf{Symbol} & \textbf{Assumptions} & \textbf{Dim.} & \textbf{Parameters (Thm.~\ref{ltp-3})} \\
 \hline

	& $\displaystyle\varrho(x) \frac{|\xi|^2 (|\xi|^2+1)}{(|\xi|^2+2) \log(|\xi|^2+2)}$
	& $\varrho \in C^{>0}((0,\infty))$
	& $d \geq 1$
	& \vtop{\hbox{\strut $\gamma_0(\varrho(x))=\gamma_{\infty}(\varrho(x))=2$} \hbox{\strut $m \in (0,\sqrt{2})$}} \B\T \\ \hline

	& $\displaystyle|\xi| \arctan(\varrho(x) |\xi|)$
	& $\varrho \in C^{>0}((0,\infty))$
	& $d \geq 1$
	& $\gamma_0(\varrho(x))=\gamma_{\infty}(\varrho(x))=1$, $m=0$ \B\T \\ \hline

	  TLP-like
	& $\displaystyle(|\xi|^2+\varrho(x)^2)^{\alpha(x)/2} \cos \big(\alpha(x) \arctan \tfrac{|\xi|}{\varrho(x)}\big)-\varrho(x)^{\alpha(x)}$
	& \vtop{\hbox{\strut $\alpha \in C^{>0}((0,2))$} \hbox{\strut $\varrho \in C^{>0}((0,\infty))$}}
	& $d \geq 1$
	& \vtop{\hbox{\strut $\gamma_0(\alpha(x),\varrho(x))=2$, $m \in (0,\varrho^L)$} \hbox{\strut $\gamma_{\infty}(\alpha(x),\varrho(x))=\alpha(x)$}}  \B\T \\ \hline

	& $\displaystyle \varrho(x) |\xi| \frac{\cosh^2(\sqrt{2} |\xi|)}{\sinh(2 \sqrt{2} |\xi|)}$
	& $\varrho \in C^{>0}((0,\infty))$
	& $d \geq 1$
	& \vtop{\hbox{\strut $\gamma_0(\varrho(x))=\gamma_{\infty}(\varrho(x))=1$} \hbox{\strut $m=0$}}  \B\T \\ \hline

	& $\displaystyle \varrho(x) |\xi| \frac{\sinh^2(\sqrt{2} |\xi|)}{\sinh(2 \sqrt{2} |\xi|)}$
	& $\varrho \in C^{>0}((0,\infty))$
	& $d \geq 1$
	& \vtop{\hbox{\strut $\gamma_0(\varrho(x))=\gamma_{\infty}(\varrho(x))=1$} \hbox{\strut $m=0$}} \B\T \\ \hline

	& $\varrho(x) |\xi| \coth((2 |\xi|)^{-1})-\varrho(x) |\xi|^2$
	& $\varrho \in C^{>0}((0,\infty))$
	& $d \geq 1$
	& $\gamma_0(\varrho(x))=1$, $\gamma_{\infty}(\varrho(x))=2$, $m=0$ \B\T \\ \hline

	& $\displaystyle\varrho(x) \log(\sinh(\sqrt{2} |\xi|))-\varrho(x) \log(\sqrt{2} |\xi|)$
	& $\varrho \in C^{>0}((0,\infty))$
	& $d \geq 1$
	& $\gamma_0(\varrho(x))=2$, $\gamma_{\infty}(\varrho(x))=1$, $m=0$ \B\T \\ \hline

	 isotropic Meixner-like
	& $\displaystyle\varrho(x) \log(\cosh(\sqrt{2} |\xi|))$
	& $\varrho \in C^{>0}((0,\infty))$
	& $d \geq 1$
	&  $\gamma_0(\varrho(x))=2$, $\gamma_{\infty}(\varrho(x))=1$, $m=0$ \B\T \\ \hline

	& $\displaystyle|\xi| \log(1+\varrho(x) \tanh(b(x) |\xi|))$
	& \vtop{\hbox{\strut $b \in C^{>0}((0,\infty))$} \hbox{\strut $\varrho \in C^{>0}((0,\infty))$}}
	& $d \geq 1$
	& \vtop{\hbox{\strut $\gamma_0(b(x),\varrho(x))=\gamma_{\infty}(b(x),\varrho(x))=1$} \hbox{\strut $m=0$}} \B\T \\ \hline

	& $\displaystyle\frac{\Gamma(\varrho(x) |\xi|^2+1/2)}{\Gamma(\varrho(x) |\xi|^2)}$
	& $\varrho \in C^{>0}((0,\infty))$
	& $d \geq 1$
	& \vtop{\hbox{\strut $\gamma_0(\varrho(x))=2$, $\gamma_{\infty}(\varrho(x))=1$} \hbox{\strut $m \in (0,1/\sqrt{2\varrho^U})$}} \B\T \\ \hline

	& $\displaystyle|\xi|^2 \frac{\Gamma(\alpha(x) |\xi|^2+1-\alpha(x))}{\Gamma(\alpha(x) |\xi|^2+1)}$
	& $\alpha \in C^{>0}((0,1))$
	& $d \geq 1$
	& \vtop{\hbox{\strut $\gamma_0(\alpha(x))=2$, $\gamma_{\infty}(\alpha(x))=2-\alpha(x)$} \hbox{\strut $m \in (0,\sqrt{1/\alpha^U-1})$}} \B\T \\ \hline

	& $\displaystyle\frac{\Gamma(\alpha(x) |\xi|^2+1)}{\Gamma(\alpha(x) |\xi|^2+1-\alpha(x))} - \frac{1}{\Gamma(1-\alpha(x))}$
	& $\alpha \in C^{>0}((0,1))$
	& $d \geq 1$
	& \vtop{\hbox{\strut $\gamma_0(\alpha(x))=2$, $\gamma_{\infty}(\alpha(x))=\alpha(x)$} \hbox{\strut $m \in (0,\sqrt{1/\alpha^U})$}} \B\T \\ \hline

	 Lamperti stable-like 
	& $\displaystyle\frac{\Gamma(|\xi|^2+\alpha(x)+\varrho(x))}{\Gamma(|\xi|^2+\varrho(x))} - \frac{\Gamma(\alpha(x)+\varrho(x))}{\Gamma(\varrho(x))}$
	& \vtop{\hbox{\strut $\alpha \in C^{>0}((0,1))$} \hbox{\strut $\varrho \in C^{>0}((0,\infty))$}}
	& $d \geq 1$
	& \vtop{\hbox{\strut $\gamma_0(\alpha(x),\varrho(x))=2$,$m \in (0,\sqrt{\varrho^L+\alpha^L})$} \hbox{\strut $\gamma_{\infty}(\alpha(x),\varrho(x))=2\alpha(x)$}} 
\end{tabular} \caption{(cont.)}
\end{table}
\end{landscape}

Theorem~\ref{ltp-7} allows us to deduce an uniqueness and existence result for solutions of L\'evy-driven stochastic differential equations, i.\,.e. SDEs of the form \begin{equation*}
	dX_t = b(X_{t-}) \, dt + \sigma(X_{t-}) \, dL_t, \qquad X_0 = x \in \mathbb{R}^d
\end{equation*}
where $(L_t)_{t \geq 0}$ is a $n$-dimensional L\'evy process. \emph{If} the SDE has a unique weak solution, then it is possible to give conditions in terms of $\sigma$ and the L\'evy measure $\nu$ of $(L_t)_{t \geq 0}$ which ensure that the solution is a rich Feller process, cf.\ K\"{u}hn \cite{sde}. It is, however, in general a non-trivial problem to prove the uniqueness of the solution, see  \cite{stroock} for the particular case that $(L_t)_{t \geq 0}$ has a non-vanishing diffusion part and e.\,g.\ \cite{kul15-2,zan} for the case that $(L_t)_{t \geq 0}$ is an isotropic stable process. \par
Using the parametrix construction, we can give sufficient conditions in terms of the characteristic exponent $\psi$ such that the SDE has a unique weak solution which is a rich Feller process.
	
\begin{kor} \label{ltp-15}
	Let $(L_t)_{t \geq 0}$ be a one-dimensional L\'evy process with characteristic exponent $\psi$. Suppose that $\psi$ has a holomorphic extension $\Psi$ to $\Omega=\Omega(m^L,\theta)$ for some $m^L \geq 0$, $\theta \in (0,\pi/2)$ which satisfies the following two growth conditions: \begin{enumerate}
		\item There exist $\alpha \in (0,2]$, $\beta  \in (1,2)$ and constants $c_1,c_2>0$ such that \begin{equation*}
				\re \Psi(z) \geq c_1 |\re z|^{\beta} \fa |z| \gg 1, \; z \in \Omega
		\end{equation*}
		and \begin{equation*}
			|\Psi(z)| \leq c_2 (|z|^{\alpha} \I_{\{|z| \leq 1\}} + |z|^{\beta} \I_{\{|z|>1\}}), \qquad z \in \Omega.
		\end{equation*}
		\item There exists a constant $c_3>0$ such that $|\Psi'(z)| \leq c_3 |z|^{\beta-1}$ for all $z \in \Omega$, $|z| \gg 1$.
	\end{enumerate}
	Let $b: \mbb{R} \to \mbb{R}$ and $\sigma: \mbb{R} \to \mbb{R}$ be H\"{o}lder continuous bounded functions such that \begin{equation*}
		0 < \sigma^L := \inf_{x \in \mbb{R}} |\sigma(x)| \leq \sup_{x \in \mbb{R}} |\sigma(x)| =: \sigma^U < \infty.
	\end{equation*}
	Then there exists a unique weak solution to the SDE \begin{equation}
		dX_t = b(X_{t-}) \, dt + \sigma(X_{t-}) \, dL_t, \qquad X_0 = x, \label{sde}
	\end{equation}
	and the solution is a rich L\'evy-type process with symbol $q(x,\xi) = -ib(x) \xi + \psi(\sigma(x) \xi)$. The solution has the following additional properties: \begin{enumerate}
		\item The transition probability $p:(0,\infty) \times \mbb{R} \times \mbb{R} \to [0,\infty)$ is continuous, differentiable with respect to $t$ and satisfies the heat kernel estimates from Theorem~\ref{ltp-3} with $\gamma_0(b(x),\sigma(x)) = \min\{\alpha,1\}$, $\gamma_{\infty}(b(x),\sigma(x)) = \beta$ and any $m\in (0,m^L/\sigma^U)$. 
		\item $C_c^{\infty}(\mbb{R})$ is a core for the generator $(L,\mc{D}(L))$ of $(X_t)_{t \geq 0}$ and $C_{\infty}^2(\mbb{R}) \subseteq \mc{D}(L)$. Moreover, $p$ is a fundamental solution to the Cauchy problem for the operator $\partial_t - L$.
		\item $(X_t)_{t \geq 0}$ is the unique solution to the $(L,C_c^{\infty}(\mbb{R}^d))$-martingale problem.
		\item The associated semigroup has the strong Feller property.
\end{enumerate} \end{kor}

Corollary~\ref{ltp-15} applies, in particular, to Lévy processes $(L_t)_{t \geq 0}$ with the following characteristic exponents:  \begin{enumerate}
	\item (isotropic stable) $\psi(\xi) = |\xi|^{\alpha}$, $\xi \in \mbb{R}$, $\alpha \in (1,2]$,
	\item (relativistic stable) $\psi(\xi) = (|\xi|^2+\varrho^2)^{\alpha/2}-\varrho^{\alpha}$, $\xi \in \mbb{R}$, $\varrho>0$, $ \alpha \in (1,2)$,
	\item (Lamperti stable) $\psi(\xi) = (|\xi|^2+\varrho)_{\alpha}-(\varrho)_{\alpha}$, $\xi \in \mbb{R}$, $\varrho>0$, $\alpha \in (1/2,1)$, where $(r)_{\alpha} := \Gamma(r+\alpha)/\Gamma(r)$ denotes the Pochhammer symbol,
	\item (truncated L\'evy process) $\psi(\xi) = (|\xi|^2+\varrho^2)^{\alpha/2} \cos(\alpha \arctan(\varrho^{-1} |\xi|))-\varrho^{\alpha}$, $\xi \in \mbb{R}$, $\alpha \in (1,2)$, $\varrho>0$,
	\item (normal tempered stable) $\psi(\xi) = (\kappa^2+(\xi-ib)^2)^{\alpha/2}-(\kappa^2-b^2)^{\alpha/2}$, $\xi \in \mbb{R}$, $\alpha \in (1,2)$, $b>0$, $|\kappa|>|b|$.
\end{enumerate}
Up to know, this result was only known for the particular case that $(L_t)_{t \geq 0}$ is an isotropic stable process, see Knopova \& Kulik \cite{kul15,kul15-2} and the references therein. \par \medskip

We close this section with an existence result for L\'evy-type processes with symbols of variable order.

\begin{thm} \label{ltp-17}
	Let $I \subseteq \mbb{R}^n$ an open convex set and $\psi_{\beta}: \mbb{R}^d \to \mbb{C}$, $\beta \in I$, be a family of continuous negative definite functions satisfying \eqref{LTP1}-\eqref{LTP3} on  \begin{equation*}
				\Omega(\vartheta) := \{z \in \mbb{C} \backslash \{0\};  \arg z \in (-\vartheta,\vartheta) \cup (\pi-\vartheta,\pi+\vartheta)\}
			\end{equation*}
		for some $\vartheta \in (0,\pi/2)$. Assume, in addition, that \begin{description}
		\item[\normalfont{(LTP4')}]\label{LTP4'} The partial derivative $\frac{\partial}{\partial \beta_j} \Psi_{\beta}(r)$ exists for all $r \in \mbb{R}$ and extends holomorphically to $\Omega(\vartheta)$ for all $j \in \{1,\ldots,n\}$ and $\alpha \in I$. There exist an increasing slowly varying (at $\infty$) function $\ell: (0,\infty) \to (0,\infty)$ and a constant $c_4>0$ such that \begin{equation*}
			\left| \frac{\partial_{\beta_j} \Psi_{\beta}(z)}{\Psi_{\beta}(z)} \right| \leq c_4 (1+\ell(|z|)) \fa z \in \Omega(\vartheta), \; j=1,\ldots,n,
		\end{equation*}
	\end{description}
	and 
	\begin{enumerate}[label*=\upshape (S),ref=\upshape S]
		\item\label{S} $(\psi_{\beta})_{\beta \in I}$ satisfies the sector condition, i.\,e.\ there exists a constant $c>0$ such that  \begin{equation*}
			|\im \Psi_{\beta}(z)| \leq c |\re \Psi_{\beta}(z)| \fa z \in \Omega(\vartheta), \; \beta \in I.
		\end{equation*}
	\end{enumerate}
	Then for any two H\"{o}lder continuous mappings $\alpha:\mbb{R}^d \to (0,1]$ and $\beta: \mbb{R}^d \to I$ such that $\alpha^L := \inf_{x \in \mbb{R}^d} \alpha(x)>0$, there exists a rich L\'evy-type process $(X_t)_{t \geq 0}$ with symbol \begin{equation*}
		q(x,\xi) := (\psi_{\beta(x)}(\xi))^{\alpha(x)}, \qquad x,\xi \in \mbb{R}^d.
	\end{equation*}
	The process $(X_t)_{t \geq 0}$ has the following properties: \begin{enumerate}
		\item The transition probability $p:(0,\infty) \times \mbb{R}^d \times \mbb{R}^d \to [0,\infty)$ is continuous, differentiable with respect to $t$ and satisfies the heat kernel estimates from Theorem~\ref{ltp-5} with $\tilde{\gamma}_0(\alpha(x),\beta(x)) := \alpha(x) \gamma_{0}(\beta(x))$, $\tilde{\gamma}_{\infty}(\alpha(x),\beta(x)) := \alpha(x) \gamma_{\infty}(\beta(x)) $ and $m=0$; here $\gamma_0(\beta(x))$ and $\gamma_{\infty}(\beta(x))$ are the mappings associated with $(\psi_{\beta})_{\beta \in I}$ by the growth condition \eqref{LTP3}.
		\item $C_c^{\infty}(\mbb{R}^d)$ is a core for the generator $(L,\mc{D}(L))$ of $(X_t)_{t \geq 0}$ and $C_{\infty}^2(\mbb{R}^d) \subseteq \mc{D}(L)$. Moreover, $p$ is a fundamental solution to the Cauchy problem for the operator $\partial_t - L$.
		\item $(X_t)_{t \geq 0}$ is the unique solution to the $(L,C_c^{\infty}(\mbb{R}^d))$-martingale problem.
		\item The associated semigroup has the strong Feller property.
		\end{enumerate}
	\end{thm}

In dimension $d=1$ Theorem~\ref{ltp-17} remains valid if we just assume that $(\psi_{\beta})_{\beta \in I}$ satisfies \eqref{LTP2}, \eqref{LTP3}, \eqref{LTP4'} and \eqref{S}, i.\,e.\ we can drop the assumption of rotational invariance.

\begin{ack}
	I would like to thank Ren\'e Schilling for helpful comments and suggestions. 
\end{ack}

\end{document}